\documentclass[12pt]{amsart}
\usepackage{amssymb}
\usepackage{epsf}

\newtheorem{thm}{Theorem}[section]
\newtheorem{prop}[thm]{Proposition}
\newtheorem{lem}[thm]{Lemma}

\newtheorem{definition}[thm]{Definition}
\newenvironment{defn}{\begin{definition}\sl}{\end{definition}}
\newtheorem{remark}[thm]{Remark}
\newenvironment{rem}{\begin{remark}\rm}{\end{remark}}
\newtheorem{example}[thm]{Example}
\newenvironment{ex}{\begin{example}\rm}{\end{example}}

\numberwithin{equation}{section}

\newcommand{\subdot}{{\,\subset\!\!\!\!\cdot\,\,}}
\newcommand{\subdoteq}{\subseteq \kern -.6em\raise.1em\hbox{$\cdot$}\kern .5em}

\def\urladdr#1{}

\newcommand{\QED}{
\setlength{\unitlength}{1.0pt}%
\begin{picture}(7.5,7.5)
\put(0,-5){\rule{2.5pt}{2.5pt}}
\put(0,-2.5){\rule{5pt}{2.5pt}}
\put(0,0){\rule{5pt}{2.5pt}}
\put(2.5,2.5){\rule{5pt}{2.5pt}}
\end{picture}\vspace{10pt}}

\def\L{{\mathcal{L}}}
\newcommand{\IR}{{\mathcal{R}}}
\newcommand{\IC}{{\mathcal{C}}}
\newcommand{\HO}{{\mathcal{H}}}
\newcommand{\J}{{\mathcal{J}}}
\newcommand{\RS}{{\mathcal{S}}}
\newcommand{\RE}{{\mathcal{E}}}
\newcommand{\CE}{{\mathcal{X}}}

\newcommand{\wbar}[1]{{\overline{#1}}}

\begin{document}

\title[Complementary Algorithms for Tableaux]{Complementary Algorithms 
       for Tableaux}  
\author{Tom Roby} 
\address[Tom Roby]{Department of Mathematics\\
	California State University\\
 	Hayward, California 94542-3092\\
	USA. \tt troby@csuhayward.edu}
\urladdr{http://seki.csuhayward.edu}
\author{Frank Sottile}
\address[Frank Sottile]{Department of Mathematics and Statistics\\
        University of Massachusetts\\
        Amherst, Massachusetts  01002\\
        USA. \tt sottile@math.umass.edu}
\urladdr{http://www.math.umass.edu/\~{}sottile}
\author{Jeff Stroomer}
\address[Jeff Stroomer]{Xilinx, Incorporated\\
         2300 55th Street\\
         Boulder, Colorado 80301\\
         USA. \tt Jeff.Stroomer@xilinx.com }
\author{Julian West}
\address[Julian West]{Department of Mathematics and Statistics\\
	University of Victoria\\
	Victoria, BC V8W 3P4, Canada. \tt westj@mala.bc.ca}
\urladdr{http://www.mala.bc.ca/\~{}westj}

\thanks{25 November 2000}
\thanks{Second author supported in part by NSERC grant OGP0170279 and NSF
        grant DMS-9701755}
\thanks{Fourth author supported in part by NSERC grant OGP0105492}
\thanks{Extended version of ``Jeux de Tableaux''~\cite{fpsac00}, 
presented at FPSAC Conference, Moscow, 2000.} 
\thanks{J.~Combin.~Th., Ser.~A, to appear}

\begin{abstract}
We study four operations defined on pairs of tableaux.  
Algorithms for the first three involve the familiar procedures
of jeu de taquin, row insertion, and column insertion.
The fourth operation, hopscotch, is new, although
specialised versions have appeared previously. 
Like the other three operations, this new operation may be computed with
a set of local rules in a growth diagram, 
and it preserves Knuth equivalence class.  
Each of these four
operations gives rise to an {\em a priori} distinct
theory of dual equivalence.  
We show that these four theories coincide.  
The four operations are linked via the involutive tableau operations of
complementation and conjugation.
\end{abstract}

\maketitle

\section*{Introduction}
Sch\"utzenberger's theory of jeu de taquin~\cite{Schu63} and the insertion
procedure of Schensted~\cite{Schensted} have found their way into the 
standard toolkit of the combinatorial representation theorist.  
These algorithms were originally developed for tableaux of partition
shape, but more recent work~\cite{SS90, Ha92} uses them to define
operations on pairs of skew tableaux.
A striking duality between Schensted insertion and Sch\"utzenberger's jeu de
taquin was noted by Stembridge~\cite{Stem87} in 
his theory of rational tableaux.

We introduce the involutive tableau operation of complementation to
formalise this duality, showing that internal
row insertion and the jeu de taquin give complementary operations on pairs
of tableaux.
This is readily seen using the growth diagram formulation for operations on
tableaux introduced by Fomin~\cite{Fomin_86} and rediscovered several times,
most notably by van Leeuwen~\cite{Leeuwen_CWI,Leeuwen_EJC-3}.
Column insertion~\cite{Fu97,Sagan} is another tableau algorithm that
can be extended 
to define the operation of internal column insertion between pairs of
tableaux. 
The operation hopscotch which is complementary to internal column insertion is
new, although it extends both Stroomer's algorithm of 
column-sliding~\cite{Stroomer} and
Tesler's~\cite{Tesler} rightward- and leftward-shift games. 

These four operations of jeu de taquin, internal row insertion,  internal
column insertion, and hopscotch all preserve Knuth-equivalence.
A consequence is that hopscotch gives rise to new algorithms to rectify
tableaux. 
We also show that these four operations have the same dual-equivalence theory.

This paper is organised as follows.
Section~\ref{Tableaux}  recalls basic definitions and terminology concerning
tableaux, the jeu de taquin, and Knuth equivalence.
In Section~\ref{Complementation} we introduce complementation and determine
its behaviour with respect to Knuth equivalence and dual equivalence.
Section~\ref{LR} is devoted to operations on pairs of tableaux given by 
growth diagrams constructed from local rules, 
and introduces the notion of complementary operations on pairs of tableaux.
In Section \ref{IR} we show that the
two operations of jeu de taquin and internal row insertion are 
complementary.  
In Section \ref{GNU} we define internal column
insertion and its complementary operation hopscotch, and in Section 6, we
investigate some features of hopscotch.
Finally, in Section 7, we supplement the jeu de taquin, describing four
additional algorithms to rectify skew tableaux.
The first three are not well-known, while the fourth is new.

Schematically, the four operations of jeu de taquin, internal row
insertion, internal column insertion, and hopscotch are linked through
complementation and conjugation as follows:
$$
\begin{picture}(390,94)(-65,-15)
\thicklines
\put( 0, 69){Jeu de taquin}  \put(203, 69){Hopscotch}
\put(11,  0){Internal}       \put(210,  0){Internal}
\put(-5,-12){row insertion}  \put(188,-12){column insertion}
\put( 31,15){\vector( 0, 1){48}}
\put( 31,45){\vector( 0,-1){33}} \put(-64,35){Complementation}
\put(230,15){\vector( 0, 1){48}}
\put(230,45){\vector( 0,-1){33}} \put(235,35){Complementation}
\put( 85,-3){\vector( 1, 0){95}}
\put(105,-3){\vector(-1, 0){30}} \put(95,4){Conjugation}
\end{picture}
$$

\section{Tableaux}\label{Tableaux}

The tableaux operations of internal row insertion,
internal column insertion, and hopscotch require us to extend the usual notion
of tableaux. 
While these operations are defined on a class of more general
tableaux, we show in Section 7 how they provide new algorithms for ordinary
tableaux. 
Rather than give a common generalisation, we instead give
a definition that will suffice until 
Section~\ref{GNU}, where we will be precise about necessary further
extensions to our notion of tableau. 

A partition $\alpha$ is a weakly decreasing sequence of
positive integers $\alpha_m \geq \cdots \geq \alpha_n$ indexed by an
interval $m, m{+}1,\ldots,n$ of integers.
The shape of a partition $\alpha$ is a left justified array of
boxes whose $i$th row contains $\alpha_i$ boxes.  
We make no distinction between a partition and its shape. 
Thus we identify partitions that differ only in their number of trailing
zeroes. 
Partitions are partially ordered $\subseteq$ 
by componentwise comparison of sequences.
If $\alpha\subseteq\beta$, then we may form the {\em (skew) shape}
$\beta/\alpha$, consisting 
of those boxes which are in $\beta$ but not in $\alpha$.  
The shape $\beta/\alpha$ has {\it inner border} $\alpha$ and {\it outer
border} $\beta$.
We call $\beta/\alpha$ a {\it horizontal strip} if no column of
$\beta/\alpha$ contains two or more boxes.
That is,
$\beta_m\geq\alpha_m\geq\beta_{m+1}\geq\cdots\geq\beta_n\geq\alpha_n$.
When $\alpha\subseteq\beta\subseteq\delta$, we say that the shape 
$\delta/\beta$ {\it extends} the shape $\beta/\alpha$.

A {\it tableau} $T$ with shape $\beta/\alpha$ and entries from the
alphabet $[k]:=\{1,\ldots,k\}$ is a chain 
$\alpha =\beta^{0}\subseteq \beta ^{1}\subseteq \cdots\subseteq
\beta^{k}=\beta$,  where the
successive skew shapes $\beta^{i}/\beta^{i{-}1}$ are horizontal strips.
The {\it content} of $T$ is the sequence whose $i$th component is the number
of boxes in the horizontal strip
$\beta^{i}/\beta^{i{-}1}$. 
Filling the boxes of $\beta^{i}/\beta^{i{-}1}$ with the integer
$i$ shows this is equivalent to the usual definition of a column-strict
tableau. 
We represent tableaux both as chains and as fillings of a shape with a
totally ordered alphabet (in practice, the positive integers or Roman
alphabet).  
For example, the chain and the filling of the shape $(4,2,1)/(1)$ below
both represent the same tableau.
$$
  \epsfxsize=325.33pt \epsfbox{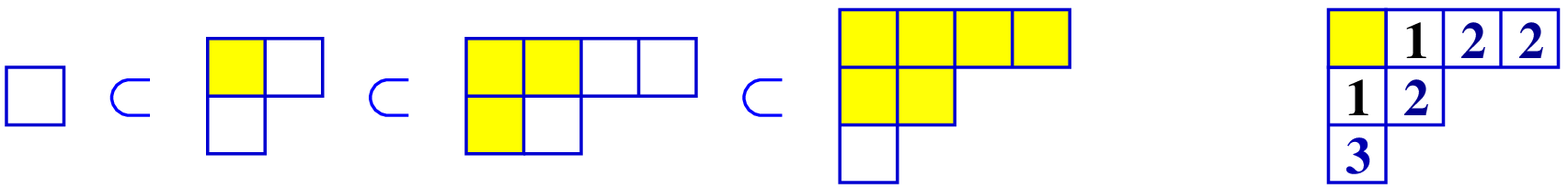}
$$

A tableau is {\it standard} if each $\beta^{i}/\beta^{i{-}1}$ consists of a
single box, so that the chain is saturated.  In this case, we say
$\beta^{i{-}1}\subdot\beta^{i}$ is a {\em cover}.   
Many (but not all) of our tableaux algorithms may be computed via the 
{\it standard renumbering} of a tableau $T$.
The standard renumbering of $T$ is the refinement of the chain representing
$T$ where each horizontal strip is filled in `left-to-right'.
For example, here is the standard renumbering of the tableau above.
$$
  \epsfxsize=49.33pt \epsfbox{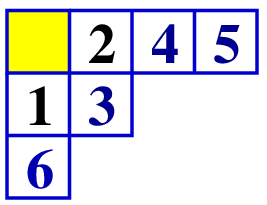}
$$

Let $\alpha\subseteq\beta\subseteq\delta$ be partitions.
If $\alpha\subdot\beta$ is a cover, then $\beta/\alpha$ is a single
box $b$, 
which we call an {\em inner corner} of $\delta/\beta$.  
If instead $\beta\subdot\delta$, then $\delta/\beta$ is a single box
$b$, called an {\em outer corner} of $\beta/\alpha$.

A {\it jeu de taquin slide}~\cite{Schu63} is a specific reversible procedure
that, given a tableau $T$ and an inner corner $b$ of the shape of $T$,
moves $b$ through $T$, producing a new tableau 
together with an outer corner. 
It proceeds by successively interchanging the empty box $b$ with one of its
neighbours to the right or below in such a way that at every step the figure
has increasing columns and weakly increasing rows.
When the box has neighbours both right and below it moves as indicated:
\begin{equation}\label{eq:js-def}
 \raisebox{-10pt}{\epsfxsize=25.33pt\epsfbox{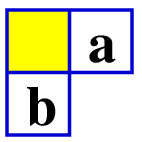}}\quad
 \longmapsto \quad \left\{\begin{array}{lcl}
  \raisebox{-10pt}{\epsfxsize=25.33pt\epsfbox{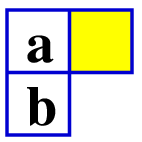}}&\quad&a<b\\
  \raisebox{-10pt}{\epsfxsize=25.33pt\epsfbox{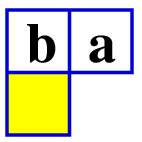}}&&b\leq  a
  \rule{0pt}{25pt}
  \end{array}\right.
\end{equation}
When the box has only a single such neighbour, it interchanges with that
neighbour, and the slide concludes when the box has no neighbours.
The reverse of this procedure is also called a jeu de taquin slide.

Sch\"utzenberger's jeu de taquin is the procedure that, given a
tableau $P$, applies jeu de taquin slides beginning at inner
corners of $P$, and concludes when there are no more such corners.
Sch\"utzenberger proves~\cite{Schu77} that the result, which we call the
{\it rectification} of $P$, 
is independent of choices of inner corners.

A tableau $U$ {\it extends} another
tableau $T$ if the shape of $U$ extends the shape of $T$.  
If $T$ is a tableau extended by $U$, then the standard renumbering
of $T$ gives a set of instructions for applying jeu de taquin slides to
$U$:
The last cover of $T$ is an inner corner of $U$, and after applying the
slide beginning with that corner, the next slide
begins at the new inner corner given by the next to last cover of $T$, and
so on.
Write $Q$ for the resulting tableau.
One could consider the sequence of outer corners obtained from
this procedure as a tableau $P$ or regard $U$ as a set of 
instructions for slides on $T$; 
these both produce the same tableau~\cite{Ha92,BSS} (up to standard
renumbering).  
If we set $\J(T,U)=(P,Q)$, then $\J$ defines an
involution, which we call
the {\it jeu de taquin}, on pairs of tableaux
where one extends the other. 

Two fundamental equivalence relations among tableaux are Knuth equivalence
and dual equivalence.
We call two tableaux $T$ and $U$ {\it Knuth-equivalent}\/ if one can
be obtained from the other by a sequence of jeu de taquin slides.
(This is equivalent to their reading words being Knuth-equivalent
in the standard sense~\cite{Stanley_ECII}.)
Two tableaux $T$ and $U$ with the same shape are {\it dual equivalent}\/ if
applying the same sequence of jeu de taquin slides to $T$ and to $U$
always gives tableaux of the same shape. 
These are related by a result of Haiman~\cite{Ha92}:
The intersection of any Knuth-equivalence class and any dual-equivalence
class is either empty or consists of a unique tableau.

\section{Complementation}\label{Complementation}

Complementation originated in a combinatorial procedure of Stanley~\cite{St83}
to model the evaluation 
$(x_1\cdots x_k)^{\lambda_1}s_\lambda(1/x_1,\ldots,1/x_k)$, where 
$s_\lambda(x_1,\ldots,x_k)$ is the Schur polynomial~\cite{Stanley_ECII}.
In combinatorial representation theory, complementation provides a
combinatorial model for the procedure of dividing by the
determinantal representation of $GL_k$.
In this context, it was used by Stembridge~\cite{Stem87} and
Stroomer~\cite{Stroomer} to develop combinatorial algorithms for studying 
{\sl rational} representations of $GL_k$.
(The classical  Robinson-Schensted-Knuth
correspondence~\cite{Ro38,Schensted,Knuth} is used to study {\sl polynomial}
representations of $GL_k$.)
Reiner and Shimozono~\cite{ReS98} studied the commutation of a version of
complementation (which they called ``Boxcomp'') with other tableaux
operations.  
Given a tableau $T$ filled with integers from $[k]$, these authors
formed a tableau whose columns were obtained from the columns of $T$ by
complementing each in the set $[k]$ and 
rotating the resulting tableau by $180^\circ$.
Rather than rotate the result, we instead choose to renumber the resulting
tableau.

More precisely, let $T$ be a tableau with shape $\beta/\alpha$ and entries
from the alphabet $[k]$, and fix $l\geq \beta_m$, the initial and largest
part of $\beta$. 
Form a new tableau $T^C$ with $l$ columns as follows.
For each $1\leq j\leq l$, let $A$ be the set-theoretic complement of the
entries of column $j$ of $T$ in the set $[k]$, considered in the dual
order: $k<k-1<\cdots <1$.  Since it is inconvenient (and possibly
ambiguous) to work with tableau whose alphabet has more than one order,
we make the replacement $j\mapsto k+1-j$, which indicates the same chain
of shapes in the usual order on $[k]$. 
The figure below illustrates this two-step process, complementing 
the tableau $T$ on the left with $k=4$ and $l=5$.
In the middle figure, we write the complement of a column of $T$ below that
column, but in reverse order, and the rightmost figure is the complement
$T^C$, where we have applied the substitution $j\mapsto k+1-j$ to the 
complemented columns and omitted writing $T$.
\begin{equation}\label{complements}
  \begin{array}{c}\epsfxsize=242.66pt \epsfbox{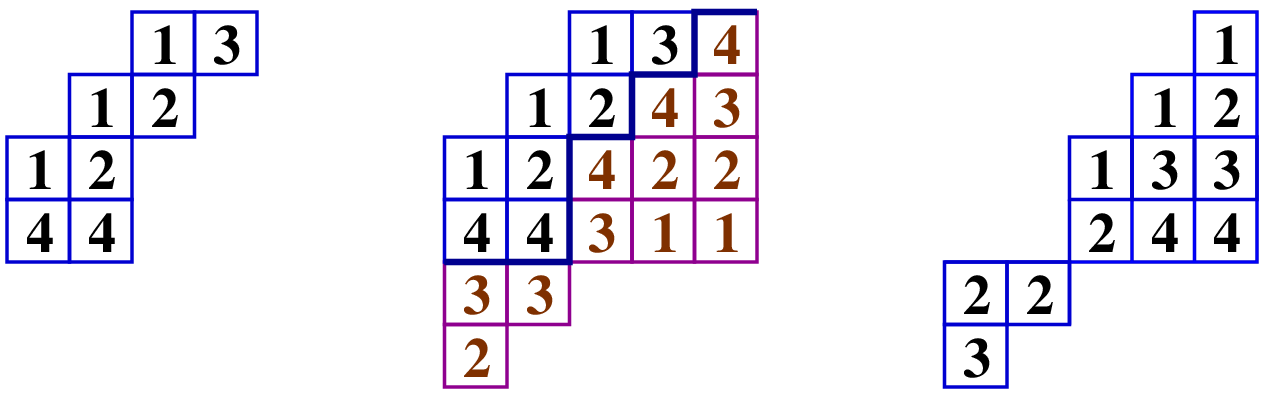}
  \end{array}
\end{equation}
This figure also illustrates the fact that if some columns of $T$ are empty 
(for instance, if $l > \beta_m$), then
the corresponding columns of $T^C$ will consist of the
full set $[k]$.
If we identify tableaux which differ by a vertical shift, we
may also write the columns of $T^C$ above the corresponding columns of $T$.
This proves that complementation is involutive, under this identification. 

\begin{thm}\label{thm_one}
Let $T$ be a tableau with shape $\beta/\alpha$ and entries from $[k]$ and
suppose $l\geq \beta_m$, the initial part of $\beta$.
Then $T^{CC}=T$.
\end{thm}

Complementation depends upon both $k$ and $l$.
Our notation, $T^C$, intentionally disregards this dependence.
We adopt the convention that we use the same integers $k$ and 
$l$ when two or more tableaux are to be complemented.
Using an extension of notation, we observe that
the shape of $T^C$ is $(l^k,\beta) / \alpha$.  (Here $l^{k}$ denotes the
rectangular shape consisting of $k$ parts of length $l$; context will
distinguish this from the superscripts on Greek letters indicating
chains of shapes.)  

\begin{figure}[htb]
$$
  \epsfxsize=350pt \epsfbox{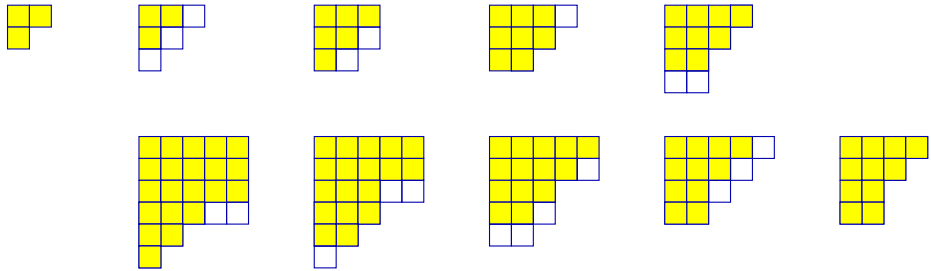}
$$
\caption{Complementary tableaux as chains.\label{fig:comp2}}
\end{figure}

We express complementation in terms of chains in Young's lattice.
Let $T$ be a tableau, written as a chain 
$\alpha = \beta^0 \subseteq 
     \beta^1 \subseteq \cdots \subseteq \beta^k = \beta$
in Young's
lattice, where $\beta^i/\beta^{i{-}1}$ is the horizontal strip of $i$'s in $T$.  
Then $T^C$ is the chain $\beta\subseteq(l,\beta^{k{-}1})\subseteq
(l^2,\beta^{k{-}2})\subseteq\cdots\subseteq 
(l^i,\beta^{k{-}i})\subseteq\cdots\subseteq(l^k,\beta^0)$.  
It is an exercise that this agrees with the definition given above.
In particular, it can immediately be seen that the shape is as claimed.
Each row of partitions in Figure~\ref{fig:comp2} is one of the two
complementary tableaux in~(\ref{complements}).
We write the second row in reverse order to illustrate the key feature
of complementation---that the $i$th horizontal strip in the complement
contains boxes in exactly the columns complementary to the corresponding
($k{+}1{-}i$)th horizontal strip in the original tableau.

Complementation preserves both Knuth equivalence and dual equivalence.

\begin{thm} \label{comp-equiv}
Suppose $T$ and $U$ are tableaux with at most $l$ 
columns and entries from the alphabet $[k]$.
Then
\begin{itemize}
\item[(i)] $U$ is Knuth-equivalent to $T$ if and only if\/
$U^C$ is Knuth-equivalent to $T^C$.
\item[(ii)] $U$ is dual-equivalent to $T$ if and only if\/
$U^C$ is dual-equivalent to $T^C$.
\end{itemize}
\end{thm}

We prove Theorem~\ref{comp-equiv} in Section~\ref{IR}. 
\medskip

\noindent{\bf Remark.}
Theorem~\ref{comp-equiv} shows that complementation commutes with the involution
reversal $T\mapsto T^e$ of~\cite{BSS}.
To see this, given a tableau $T$, let $T^*$ be the tableau obtained by rotating
$T$ $180^\circ$ and replacing each entry $j$ with $k{+}1{-}j$.
Then the {\it reversal} $T^e$ of $T$ is the unique tableau dual-equivalent
to $T$ (and hence with the same shape as $T$) and Knuth-equivalent to $T^*$.
For tableaux of partition shape, reversal coincides with Sch\"utzenberger's
evacuation procedure~\cite{Schu63} (called ``promotion'' there), but the two
procedures differ for general skew tableaux.
This extends the result of Reiner and Shimozono (\cite{ReS98}, Theorem 2) that
complementation commutes with evacuation.

\section{Growth diagrams and local rules}\label{LR}

Many properties of tableaux algorithms such as symmetry become clear when
the algorithms are formulated in terms of growth diagrams governed by local
rules.
Fomin~\cite{Fomin_86} introduced this approach to the
Robinson-Schensted correspondence, it was rediscovered by van
Leeuwen~\cite{Leeuwen_EJC-3}, and Roby~\cite{Roby} developed 
it further.  
We study tableaux algorithms related via complementation of their growth
diagrams. 

\subsection{Growth diagrams}\label{Growth}
A {\it growth diagram} is a rectangular array of
partitions where every row and every column is a tableau, with the
additional restriction that all tableaux formed by the rows have the
same content, 
as do all tableaux formed by the columns.
Specifically, the sequence $\beta^0,\beta^1,\ldots,\beta^k$ of partitions
(read left-to-right) in each row
forms a chain 
$\beta^{0}\subseteq \beta ^{1}\subseteq \cdots\subseteq \beta^{k}$,
with each $\beta^{i}/\beta^{i{-}1}$ a horizontal strip, and the number of
boxes in  $\beta^{i}/\beta^{i{-}1}$ does not depend upon which row this chain
came from. 
(We require the same to hold for the sequence of partitions read
top-to-bottom in each column.)
For example, here is a growth diagram where the horizontal tableaux have
content $(1,2,2)$ and the vertical tableaux have content $(2,1,1)$.
\begin{equation}\label{growth}
  \begin{array}{llll}
     1 & 2 & 31 & 33 \cr 
    21 &22 &321 &332 \cr 
    211&221&3211&3321\cr 
    221&222&3221&3322\cr 
  \end{array}
\end{equation}

Traditionally the tableaux in a growth diagram are standard.  We 
relax this in order to define complementation of a growth diagram.
Given an integer $l\geq\delta_1$, where $\delta$ is the lower right
partition in the growth diagram, we complement the tableaux
represented by each row to obtain new tableaux, also written as 
chains of shapes.  These combine together to give a new growth diagram.
We then complement the columns of these diagrams to obtain two more
growth diagrams.  Figure~\ref{fig:four} shows the four growth diagrams
we obtain from~(\ref{growth}) by this process.
Here, $k$ is 3 for both the vertical and horizontal tableaux, and the
complementation parameter is $l=4$.

\begin{figure}[htb]
$$
 \begin{array}{|llll|llll|}\hline
   1   & 2   & 31   & 33   &  33   & 431    & 442   & 4441  \cr 
   21  & 22  & 321  & 332  &  332  & 4321   & 4422  & 44421  \cr 
   211 & 221 & 3211 & 3321 &  3321 & 43211  & 44221 & 444211 \cr
   221 & 222 & 3221 & 3322 &  3322 & 43221  & 44222 & 444221 \cr\hline
   221 & 222 & 3221 & 3322 &  3322 & 43221  & 44222 & 444221 \cr 
   4211& 4221& 43211& 43321&  43321& 443211 & 444221& 4444211 \cr
   4421& 4422& 44321& 44332&  44332& 444321 & 444422& 4444421 \cr
   4441& 4442& 44431& 44433&  44433& 444431 & 444442& 4444441\cr\hline
 \end{array} 
$$
\caption{Complementary growth diagrams.\label{fig:four}}
\end{figure}
The lower right growth diagram may also be obtained from the lower left
diagram by complementing rows.
Indeed, if we number the rows of the original growth diagram
$1,2,\ldots,k$ and the 
columns $1,2,\ldots,m$ and if $\beta$ is the partition in position
$(i,j)$ of the 
original diagram, then $(l^{(k{-}i){+}(m{-}j)},\beta)$ is the partition in
position $(k{-}i{+}1,m{-}j{+}1)$ of the lower right growth diagram.
By Theorem~\ref{thm_one}, further complementation of the rows or columns
yields no new growth diagrams.
(Here, we identify tableaux which differ by a vertical shift.) 

\subsection{Local rules}
A local rule $\L$ is a rule for completing the missing corner of a
$2\times 2$ growth diagram given 3 partitions that are related by horizontal
strips.
Suppose we have partitions $\alpha\subset\beta\subset\delta$ with
$\beta/\alpha$ and $\delta/\gamma$ horizontal strips.
A {\it switching local rule} $\L$ is a rule for completing 
the missing lower left corner $\gamma$ of a partial growth diagram
$$
  \begin{array}{cc}\alpha&\beta\\&\delta\end{array}\,.
$$
Write 
$\gamma =\L(\alpha,\beta,\delta)$.
A switching local rule also
gives a rule for completing a missing upper right corner of a partial
growth diagram, by the obvious symmetry of the two cases. 
In order to make $\L$ reversible, we insist
that the rule be symmetrical: $\gamma =\L(\alpha,\beta,\delta) \iff
\beta =\L(\alpha ,\gamma ,\delta )$.

Suppose we have 
three partitions $\alpha,\beta$, and $\gamma$ with $\gamma/\alpha$ and
$\beta/\alpha$ horizontal strips, so that 
$$
  \begin{array}{cc}\alpha&\beta\\\gamma& \end{array}
$$
is a (partial) growth diagram.
An {\it insertion local rule} $\L$ is a rule that
associates a fourth partition $\delta$ to 
such a partial growth diagram such that 
$$
  \begin{array}{cc}\alpha&\beta\\\gamma& \delta\end{array}
$$
is a growth diagram.
We require that  $\L$ be invariant under vertical
shifts of the horizontal strips $\beta/\alpha$ and $\gamma/\alpha$.
By this we mean that 
if $l$ is at least as long as the initial
part of $\delta$, then the rule $\L$ completes the partial growth
diagram
$$
  \begin{array}{cc}(l,\alpha)&(l,\beta)\\(l,\gamma)& \end{array}
$$
with the partition $(l,\delta)$.  
We write $\delta =\L(\alpha ;\beta ,\gamma )$.  
We further require $\L$ to be symmetric in $\beta$ and  $\gamma$.

We would like to define a reverse map 
$\L(\beta ,\gamma; \delta) $ for all $\beta$, $\gamma$, and $\delta$ with 
$\delta/\beta$ and $\delta/\gamma$ horizontal strips with the property that
$\L(\alpha;\beta,\gamma)=\delta$ if and only if 
$\alpha =\L(\beta ,\gamma; \delta) $.
Unfortunately, this is impossible as the following example shows.
Let $\beta=\gamma=0$ be the empty partition
and $\delta=1$ be the partition with a single part of size 1.
Then there simply does not exist a partition $\alpha$ with $\beta/\alpha$ and
$\gamma/\alpha$ a single box.

To circumvent this problem, we call an insertion local rule $\L$ 
{\it reversible} if given partitions $\beta$, $\gamma$, and $\delta$ with 
$\delta/\beta$ and $\delta/\gamma$ horizontal strips and an integer $l$ at
least as large as the initial part of $\delta$, then there exists 
a unique partition $\alpha$ such that 
$\L(\alpha; (l,\beta),(l,\gamma))=(l,\delta)$.

A reversible local rule $\L$ and an integer $l$ together define a local rule
for computing the missing upper left corner of a growth diagram
$$
  \begin{array}{cc}&\beta\\\gamma&\delta \end{array}
$$
when $l\geq \delta_m$, the initial part of $\delta$.
First prepend a single, $(m{-}1)$th, part $l$ to each of $\beta,\gamma$, and
$\delta$ 
and then set $\L(\beta ,\gamma; \delta):=\alpha$, where 
$\alpha$ is the unique partition such that
$\L(\alpha; (l,\beta),(l,\gamma))=(l,\delta)$.
Prepending the part of size $l$ does not alter the horizontal strips
$\delta/\beta$ and $\delta/\gamma$, as it only shifts them vertically.

As with complementation, we suppress this parameter $l$ in our notation.
However, we insist that it coincides with the complementation parameter when
combining an insertion local rule with complementation.

\subsection{Tableaux algorithms}\label{TaAl}
An insertion local rule $\L$ 
determines a bijection on pairs of tableaux which share a
common border as follows.
Given a pair $(P,Q)$ sharing an inner border, write the tableau $P$ across
the first row of an 
array and the tableau $Q$ down the first column. 
Then use $\L$ to fill in the array and  obtain a
growth diagram. 
If $U$ is the tableau of the last column in this diagram and 
$T$ the tableau of
the last row in this diagram, then $T$ and $U$ share the same outer
border and the pair $(T,U)$ is determined from the
pair $(P,Q)$ by the local rule $\L$. 
When $\L$ is reversible and $P$ and $Q$ share an outer border occupying
columns (at most) $1,\ldots,l$, then we write $P$ and $Q$ 
across the last row and column of the array and use the (reverse) local rule
to complete the growth diagram, obtaining a pair $(T,U)$ which share
an inner border.  
We indicate this by writing $\L(P,Q)=(T,U)$.

Similarly, a switching local rule $\L$ 
determines an involution on pairs of Young tableaux $(P,Q)$ where 
$P$ extends $Q$ and we write $\L(P,Q)=(T,U)$.
These mappings have the following properties. 

\begin{thm}
A switching local rule $\L$ 
determines an involution 
$$
  \left\{
   \begin{minipage}[c]{1.39in}
    Tableaux $P$ and $Q$ where $Q$ extends $P$.
   \end{minipage} 
  \right\}
  \ \stackrel{\L}{\leftarrow\joinrel\relbar\joinrel\to}\ 
  \left\{
   \begin{minipage}[c]{1.39in}
    Tableaux $T$ and $U$ where $T$ extends $U$.
   \end{minipage} 
  \right\}
$$
such that if\/ $\L(P,Q)=(T,U)$, then 
$P$ and $T$ have the same content, as do $Q$ and $U$.
Also, $P$ and $U$ have the same inner border and 
$Q$ and $T$ have the same outer border.\smallskip 

A reversible insertion local rule $\L$ determines a bijection
$$
  \left\{
   \begin{minipage}[c]{1.85in}
    Tableaux $P$ and  $Q$ which share an inner border and occupy columns
   $1,\ldots,l$. 
   \end{minipage} 
  \right\}
  \ \stackrel{\L}{\leftarrow\joinrel\relbar\joinrel\to}\ 
  \left\{
   \begin{minipage}[c]{1.85in}
    Tableaux $T$ and $U$ which share an outer border and occupy columns
   $1,\ldots,l$.
   \end{minipage} 
  \right\}
$$
with $\L\circ\L$ the identity such that if\/ $\L(P,Q)=(T,U)$, then    
$P$ and $T$ have the same content, as do $Q$ and $U$.  
Also, the outer border of $P$ equals the inner border of\/ $U$, and the outer
border of $Q$ equals the inner border of\/ $T$.
\end{thm}

We combine complementation of growth diagrams
with this local rules construction of tableaux algorithms.   
Given a local rule (and mapping) $\L$, define its {\it complement}
$\L^C$ by
\begin{equation}\label{map-comp}
  \L^C(P,Q)=(T,U^{C})\mbox{, \ where \ }
  \L(P,Q^{C})=(T,U)\,,
\end{equation} 
{\sl when this is well-defined}.  
(Recall our convention that $k$ and $l$ are fixed when complementing
several tableaux.)   In other words, complement one tableau,
perform the operation determined by $\L$, then complement back
the appropiate tableau.  

One set of conditions on $\L$ which ensures this is well-defined
is the following.
\begin{itemize}
\item[I.] $\L$ does not increase the number of columns in tableaux.
    By this we mean that if $\L(P,Q)=(T,U)$ with $P$ and $Q$
    having at most $l$ columns, then $T$ and $U$ have at most $l$ columns.
    (This is automatically satisfied by switching local rules.)
\item[II.]   We have $\L(P,Q)=(T,U)$ if and only if 
         $\L(P^C,Q^C)=(T^C,U^C)$.
\end{itemize}
Any local rule $\L$ satisfying conditions I and II has a
complement $\L^C$ which also satisfies conditions I and II.

\begin{figure}[htb]
$$
\begin{picture}(185,122)(0,10)
\thicklines
\put(  0, 40){$P^C$}\put(  5, 90){$P$}
\put(173, 40){$P^C$}\put(173, 90){$P$}
\put( 77, 35){$T^C$}\put( 82, 95){$T$}
\put( 45, 75){$U$}  \put(135, 75){$U^C$}
\put( 55,  7){$Q$}  \put(125,  7){$Q^C$}
\put( 55,127){$Q$}  \put(125,127){$Q^C$}

\put(20, 20){\line(1,0){150}}   \put( 20,20){\line(0,1){100}}
\put(20, 70){\line(1,0){150}}   \put( 95,20){\line(0,1){100}}
\put(20,120){\line(1,0){150}}   \put(170,20){\line(0,1){100}}

\put( 26,24){\vector(3,2){63}}    \put( 26,116){\vector(3,-2){63}}
\put(101,74){\vector(3,2){63}}    \put(101, 66){\vector(3,-2){63}}
\put( 89,66){\vector(-3,-2){63}}  \put( 89, 74){\vector(-3,2){63}}
\put(164,116){\vector(-3,-2){63}} \put(164, 24){\vector(-3,2){63}}

\put( 47, 50){$\L^C$} \put(130, 50){$\L$}
\put( 55,100){$\L$} \put(122,100){$\L^C$}

\end{picture}
$$
\caption{Growth diagrams of complementary operations.\label{fig:4-diag}}
\end{figure}
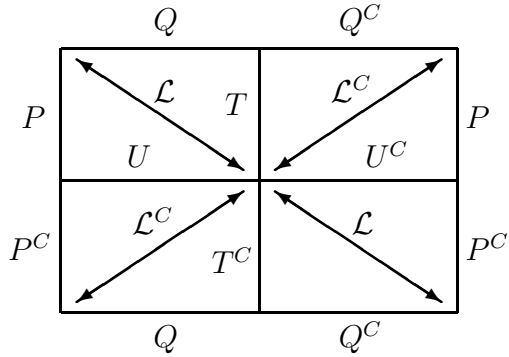

When an insertion local rule satisfies conditions I and II,  the
tableaux operation $\L^C$  is given by the switching local rule
$\L^C(\alpha ,\beta ,\delta )= \gamma $,
where $\L$ completes the partial growth diagram
$$
  \begin{array}{cc}\beta&(l,\alpha)\\\delta&\end{array}
$$
with the partition $(l,\gamma)$, for $l\geq \delta_m$, the initial part of
$\delta$. 
Also, the growth diagrams for  $\L$ and $\L^C$
fit into an array of four growth diagrams as in Section~\ref{Growth},
displayed schematically in Figure~\ref{fig:4-diag}. 
(Shown here for an insertion local rule $\L$.)

\subsection{Dual equivalence}\label{sec:dual}
Let $\L$ be a reversible insertion local rule.
Given a tableau $Q$ sharing a border (inner or outer) with a tableau $P$, set
$(T,U)=\L(P,Q)$.
We call this passage from $P$ to $T$ (which is determined by $Q$) an 
{\it $\L$-move} applied to $P$.
Two tableaux $P$ and $P'$ with the same shape are 
{\it $\L$-dual equivalent}\/ if applying the same sequence of 
$\L$-moves to $P$ and to $P'$ gives tableaux of the same shape.
When $\L$ is a switching local rule, we may
similarly define $\L$-moves and $\L$-dual equivalence.
An elementary consequence of these definitions is that if $Q$ and $Q'$ are
$\L$-dual equivalent, then the  
$\L$-moves determined by $Q$ and $Q'$ are identical.

When a local rule $\L$ has a complement $\L^C$, 
$\L$-dual equivalence coincides with 
$\L^C$-dual equivalence.

\begin{thm}\label{R-de}
Let $\L$ be a local rule.
If~(\ref{map-comp}) defines a complementary mapping $\L^C$, 
then $\L$-dual equivalence coincides with 
$\L^C$-dual equivalence.
If $\L$ satisfies condition II, then a tableau $P$ is 
$\L$-dual equivalent to
$Q$ if and only if $P^C$ is $\L$-dual equivalent to $Q^C$.
\end{thm}

\noindent{\bf Proof. }
Complementing the second coordinates of a sequence of $\L$-moves 
applied to a tableau $P$ gives a sequence of $\L^C$-moves 
applied to $P$.
Thus $\L$-dual equivalence coincides with 
$\L^C$-dual equivalence.

When $\L$  satisfies condition II, complementing both coordinates
of a sequence of  $\L$-moves applied to a tableau $P$ gives a
sequence of  $\L$-moves applied to $P^C$.
Thus $P$ is $\L$-dual equivalent to
$Q$ if and only if $P^C$ is $\L$-dual equivalent to
$Q^C$.  \QED

\section{Internal row insertion}\label{IR}

We apply the formalism of Section~\ref{LR} to show 
Sch\"utzenberger's jeu de taquin and internal row insertion (a modification
of the procedure introduced in~\cite{SS90}) are complementary operations on
pairs of tableaux.
We first define the insertion local rule $\IR$ 
for internal row insertion.
Suppose $\alpha\subset\beta$ and $\alpha\subset\gamma$ with $\beta/\alpha$
and $\gamma/\alpha$ horizontal strips, and prepending an integer if
necessary, we assume that the initial ($m$th) parts of $\alpha,\beta$, and
$\gamma$ are equal to the same number $l$.
Set $\delta_m=l$ and for $i>m$, 
\begin{equation}\label{def:iri}
  \delta_i\;\ :=\;\ \max\{\beta_i,\gamma_i\}\; + \;
                \min\{\beta_{i{-}1},\gamma_{i{-}1}\}\; -\; \alpha_{i{-}1}\,.
\end{equation}
Define $\IR(\alpha ;\beta , \gamma ):= \delta $.

\begin{lem}
The rule~(\ref{def:iri}) defines a reversible insertion local rule $\IR$.
\end{lem}

\noindent{\bf Proof. }
The rule~(\ref{def:iri}) is symmetric in $\beta$ and $\gamma$, 
and $\alpha$ is determined by
$\beta$, $\gamma$, and $\delta$, and so it is reversible.
If $\delta$ is defined by~(\ref{def:iri}), then 
$$
  \begin{array}{cc}\alpha&\beta\\\gamma&\delta\end{array}
$$
is a growth diagram.
Indeed, as $\max\{\beta_i,\gamma_i\}\leq \delta_i$, both 
$\delta/\beta$ and $\delta/\gamma$ will be horizontal strips if 
$\delta_i\leq \min\{\beta_{i{-}1},\gamma_{i{-}1}\}$.
But this follows since $\max\{\beta_i,\gamma_i\}\leq \alpha_{i{-}1}$,
as $\beta/\alpha$ and $\gamma/\alpha$ are horizontal strips.

It remains to show that 
$\sum_i\delta_i-\sum_i\gamma_i = \sum_i\beta_i-\sum_i\alpha_i$.
Let $n$ be an index such that $\alpha_n=\beta_n=\gamma_n=0$.
Since $\delta_m=l$, we have
\begin{eqnarray*}
 \sum_i\delta_i
      &=& l\;+\;\;\sum_{i=m+1}^n \max\{\beta_i,\gamma_i\}\; + \;
                \min\{\beta_{i{-}1},\gamma_{i{-}1}\}\; - \;\alpha_{i{-}1}\\
      &=& \sum_{i=m}^n \max\{\beta_i,\gamma_i\}\; + \;
                \min\{\beta_i,\gamma_i\} \;- \;\alpha_{i{-}1}\\
      &=& \sum_i\gamma_i\;+\;\sum_i\beta_i\;-\;\sum_i\alpha_i\,,
\end{eqnarray*}
which completes the proof.
\QED

We define {\it internal row insertion} to be the tableaux operation 
determined by the insertion local rule $\IR$, as in Section~\ref{TaAl}.
Observe that $\IR$ satisfies the conditions I and II of Section~\ref{TaAl},
and so it has a complement, $\IR^C$.
We show that $\IR^C$ coincides with the jeu de taquin operation 
$\J$ of Section~\ref{Tableaux}, formalising the
duality between jeu de taquin slides and row insertion discovered by
Stembridge. 
We later relate these formulations of $\J$ and $\IR$ to the
traditional descriptions found in~\cite{SS90,Stanley_ECII}.

\begin{thm}\label{J-IR}
The tableaux operation $\IR^C$  coincides with the jeu
de taquin $\J$.
\end{thm}

\begin{figure}[htb]
\begin{center}
 \begin{tabular}{ccc}
  \epsfxsize=146.67pt\epsfbox{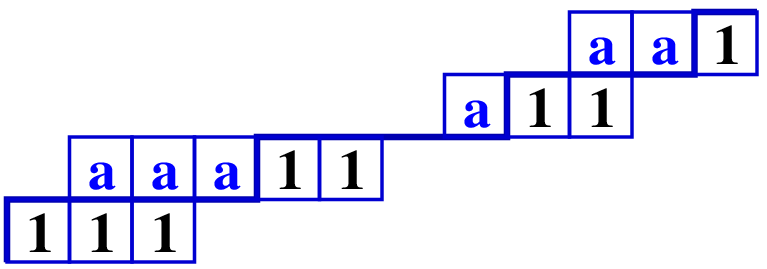}
   &\quad\raisebox{22pt}%
    {$\stackrel{\mbox{$\J$}}%
      {\leftarrow\joinrel\relbar\joinrel\longrightarrow}$}
  \quad&
  \epsfxsize=146.67pt\epsfbox{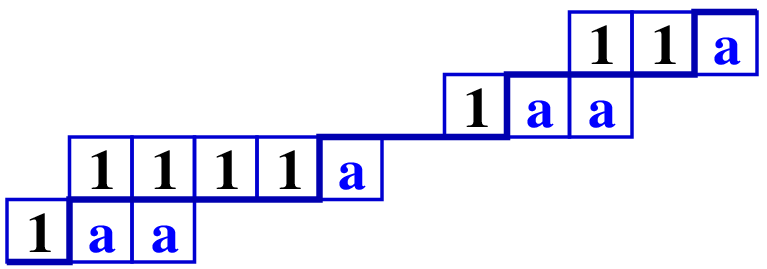}\\
  \epsfxsize=146.67pt\epsfbox{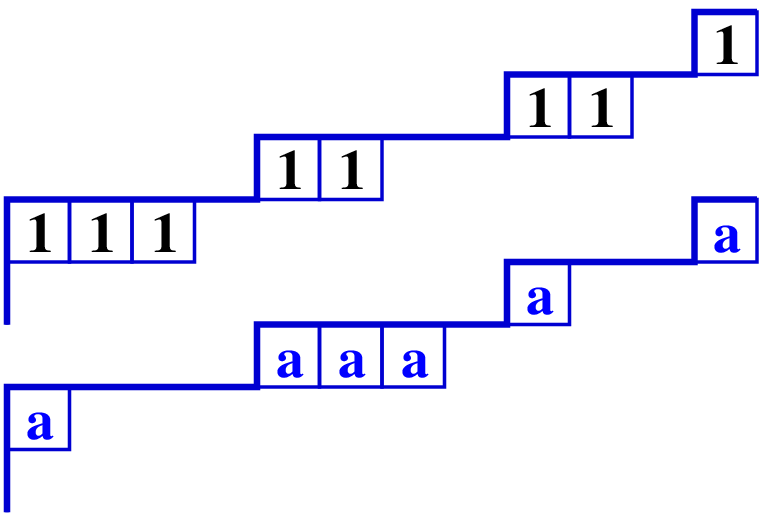}
   &\raisebox{47pt}%
    {$\stackrel{\mbox{$\IR$}}%
      {\leftarrow\joinrel\relbar\joinrel\longrightarrow}$}
  \quad&
  \epsfxsize=146.67pt\epsfbox{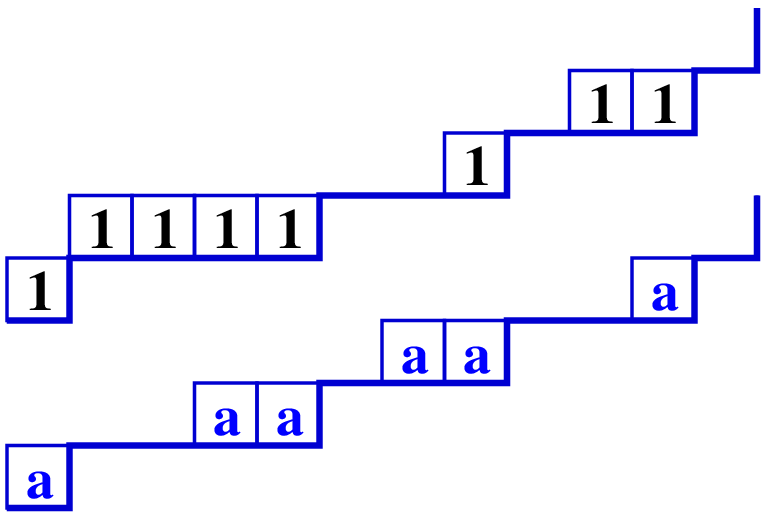}\\
 \end{tabular}
\end{center}
\caption{$\J$ and $\IR$ on complementary horizontal
strips.\label{fig:J-IRI}}
\end{figure}

Figure~\ref{fig:J-IRI} shows $\IR$ and $\J$ applied to 
complementary pairs of horizontal strips.  Our convention is to number
one tableau with letters, the other with numbers to help distinguish
them; it is only the chain of shapes that matters.  When one
tableau extends another, we juxtapose them and use a thick line to
indicate the boundary between the two.  

\medskip

\noindent{\bf Proof of Theorem~\ref{comp-equiv}. }
For statement (i), suppose we have a tableau $Q$ extending a tableau $P$.
If we complement both the rows and the columns of the growth diagram
computing  $\J(P,Q)$, then we obtain a growth diagram computing 
$\J(P^C,Q^C)$.
Thus 
$$
  (T,U)\ =\ \J(P,Q)\qquad\Longleftrightarrow\qquad
  (T^C,U^C)\ =\ \J(P^C,Q^C)\,,
$$
from which the statement (i) follows.

In the notation of this section, assertion (ii) becomes 
\begin{itemize}
\item[(ii)$'$] {\it $U$ is $\J$-dual equivalent to $T$ if and only if $U^C$
   is $\J$-dual equivalent to $T^C$.}
\end{itemize}
But this follows by Theorem~\ref{R-de}.
\QED

\noindent{\bf Proof of Theorem~\ref{J-IR}. }
Given a tableau 
$T:\alpha\subseteq\beta\subseteq\delta$ with $l=\delta_m$, the initial part
of $\delta$, we have 
$\IR^C(\alpha, \beta, \delta) = \gamma$,
where $\IR$ completes the partial growth diagram
$$
 \begin{array}{cc}\beta&(l,\alpha)\\\delta&\end{array}
$$
with the partition $(l,\gamma)$.
Thus for each $i$, 
$$
  (l,\gamma)_{i{+}1}\ =\ \max\{\delta_{i{+}1},(l,\alpha)_{i{+}1}\} + 
         \min\{\delta_i,(l,\alpha)_i\} - \beta_i
$$
or 
\begin{equation}\label{LR-jeu}
  \gamma_i\ =\ \max\{\delta_{i{+}1},\alpha_i\} + 
         \min\{\delta_i,\alpha_{i{-}1}\} - \beta_i\,.
\end{equation}

We describe this in terms of the tableau $T$.
Suppose $\beta/\alpha$ is filled with $a$'s and 
$\delta/\beta$ is filled with $1$'s.
Form a new tableau $U$ as follows.
\begin{itemize}
\item[(i)] 
     In each column of $T$ that contains both an $a$ and a $1$,
     interchange these symbols, and afterwards
\item[(ii)]
    in each row segment containing both $a$'s and $1$'s left fixed under
    (i), shift the $a$'s to the right and the $1$'s to the left.
\end{itemize}

We claim that $U$ is the tableau $\alpha\subseteq\gamma\subseteq\delta$ with
$\gamma/\alpha$ filled with $1$'s and $\delta/\gamma$ filled $a$'s.
Any tableau $U':\alpha\subseteq\gamma'\subseteq\delta$ has each column of
length 2 filled with a $1$ above an $a$, as in (i).
Observe that $\beta_i-\max\{\delta_{i{+}1},\alpha_i\}$ is the number of $a$'s
in the $i$th row of $T$ left fixed by (i) and
$\min\{\delta_i,\alpha_{i{-}1}\}-\beta_i$ the number of $1$'s left fixed by
(i). 
The claim follows by observing that $\gamma_i-\max\{\delta_{i{+}1},\alpha_i\}= 
       \min\{\delta_i,\alpha_{i{-}1}\}-\beta_i$.

The theorem follows as rules (i) and (ii) describe the action of the jeu de 
taquin on horizontal strips as given by James and
Kerber~\cite[pp.~91-92]{JK}, 
and this suffices to describe the jeu de taquin (see also~\cite{BSS}).
\QED

The local rule $\IR$ is described in terms of horizontal strips filled
with entries $a$ and $1$ as follows.
Given a horizontal strip $\beta/\alpha$ filled with $a$'s
and a horizontal strip $\gamma/\alpha$ filled with $1$'s, 
transfer any $a$'s and $1$'s which occupy the same boxes in row $i$
into the next row, beginning with $\max\{\beta_i,\gamma_i\}$.
The reverse operation is similarly described.

\begin{ex}\label{ex:iri-j}
The first column $P$ and first row $Q$ of the 
growth diagram~(\ref{growth}) are the following tableaux:
$$
  P\ =\ \raisebox{-16.5pt}{\epsfbox{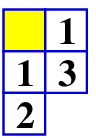}}\quad 
  {\rm and}\quad
  Q\ =\ \raisebox{-11pt}{\epsfxsize=38.67pt\epsfbox{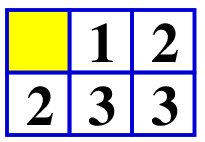}}\,.
   \smallskip
$$
We claim that this growth diagram is obtained from $P$ and $Q$ using $\IR$. 
For example, consider the upper left square of~(\ref{growth}).
Set $l:=3$, the length of the first (undisplayed) part of the partitions, and 
set  $\alpha=(3,1,0)$, $\beta=(3,2,0)$, $\gamma=(3,2,1)$, and 
$\delta=(3,2,2)$.
Note that $\delta_1=3$, and 
$$
  \begin{array}{ccccccl}
   \delta_2&=&2&=&2+3-3&=&
    \max\{\beta_2,\gamma_2\}+\min\{\beta_1,\gamma_1\}-\alpha_1\,, \\
   \delta_3&=&2&=&1+2-1&=&
    \max\{\beta_3,\gamma_3\}+\min\{\beta_2,\gamma_2\}-\alpha_2\,,
  \end{array}
$$
in agreement with~(\ref{def:iri}).
If we consider the last row and last column of~(\ref{growth}), we see that 
$\IR(P,Q)=(T,U)$, where 
$$
  T\ =\ \raisebox{-23pt}{\epsfbox{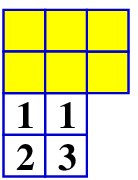}}\quad 
  {\rm and}\quad
  U\ =\ \raisebox{-23pt}{\epsfxsize=38.67pt\epsfbox{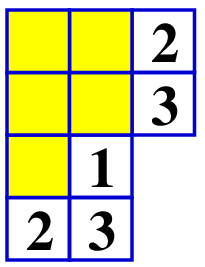}}\,.
$$

The additional 3 growth diagrams in Figure~\ref{fig:four} provide further 
examples of $\J$ and $\IR$.
Consider the lower left growth diagram of  Figure~\ref{fig:four}.
$$
  \begin{array}{llll}
   221 & 222 & 3221 & 3322  \cr 
   4211& 4221& 43211& 43321 \cr
   4421& 4422& 44321& 44332 \cr
   4441& 4442& 44431& 44433
  \end{array} 
$$
Let $P^C$ be the tableau of the first column (filled with $a<b<c$) and $Q$ the
tableau of the last row.
{}From this, we see that $\J (P^{C},Q)=(T^{C},U)$, i.e., 
$$
  \J\left(
  \raisebox{-23pt}{\epsfbox{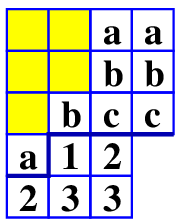}}\right)
  \quad=\quad\ 
  \raisebox{-23pt}{\epsfbox{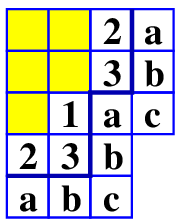}}\,.
$$
\end{ex}

Both the jeu de taquin and the Robinson-Schensted correspondence (of which
internal row insertion is a variant) are described in~\cite{Stanley_ECII} via
growth diagrams consisting of standard tableaux as follows.

In Appendix~1 to Chapter~7 in~\cite{Stanley_ECII},
Fomin gives the following local rule for the jeu
de taquin:
If $\alpha\subdot\beta\subdot\delta$, then either the interval
$[\alpha,\delta]$ in Young's lattice contains a fourth partition
$\gamma\neq\beta$, or else the interval $[\alpha,\delta]$ is a chain, in
which case we  set $\gamma:=\beta$.
One may verify from the definition~(\ref{LR-jeu}) of
$\J$ that $\gamma= \J(\alpha,\beta,\delta)$.

Similarly, suppose that $\alpha\subdot\beta$ and $\alpha\subdot\gamma$.
In Section 7.13 of~\cite{Stanley_ECII}, Stanley gives Fomin's local rule for  
Robinson-Schensted insertion.
\begin{itemize}
\item[(1)]
     If $\beta\neq \gamma$, then let $\delta:=\beta\vee\gamma$, the unique 
     partition covering both $\beta$ and $\gamma$.
\item[(2)]
     If $\beta=\gamma$, then $\beta/\alpha$ is a single box in the $i$th
     row, and we define $\delta$ so that $\delta/\beta$ is a single box in
     the $(i{+}1)$st row.
\end{itemize}
One may verify from the definition~(\ref{def:iri}) of $\IR$ that 
$\delta= \IR(\alpha;\beta,\gamma)$.
(The other 2 possibilities in~\cite{Stanley_ECII} of
$\alpha=\beta=\gamma\,(\,=\delta)$ and 
$\alpha=\beta=\gamma$ with a mark in the
square---indicating an external insertion---do not occur for us.)
\medskip

Let $st(\beta/\alpha)$ denote the standard renumbering of a horizontal
strip $\beta/\alpha$, which is now a standard tableau. 
The familiar fact that the jeu de taquin and Schensted insertion commute
with standard renumbering manifests itself here as follows:
$\delta=\IR(\alpha;\beta,\gamma)$ if and only if 
$( st(\delta/\gamma),\, st(\delta/\beta)) = 
\IR(st(\beta/\alpha),\, st(\gamma/\alpha))$, and 
$\gamma = \J(\alpha,\beta,\delta)$ if and only if 
$( st(\delta/\gamma),\, st(\gamma/\alpha)) = 
\J( st(\beta/\alpha),\, st(\delta/\beta))$.
\medskip

The skew insertion procedure of Sagan and Stanley~\cite{SS90} mixes
Schensted insertion with the internal insertion procedure $\IR$.
If tableaux $T$ and $U$ share a common inner border, then $\IR$ acts on the
pair $(T,U)$ in the same way as the forward direction of the procedure of
Theorem~6.11~\cite{SS90} when the matrix word $\pi$ is empty. 

The difference lies in the reverse procedure, where $P$ and $Q$ share an
outer border.
The essence of this difference occurs when $P$ and $Q$ are single boxes.
Suppose $\beta\subdot\delta$ and $\gamma\subdot\delta$ (so that 
$\delta/\beta$ and $\delta/\gamma$ are single boxes) and suppose we have 
$\alpha=\IR(\beta, \gamma;\delta)$.
By the definition of $\IR$, if $\beta\neq\gamma$, then $\alpha$ is the
unique partition covered by both $\beta$ and $\gamma$.
If however $\beta=\gamma$ and $\delta/\beta$ is a single box in the $i$th
row, then $\beta/\alpha$ is a single box in the ($i{-}1$)st row.

This differs from the procedure in~\cite{SS90} only in the
case (2) when $i$ is the initial row, that is, $\beta=\gamma$ and
$\delta/\beta$ is a single box in the first row. 
Then Sagan and Stanley set $\alpha=\beta$, and 
bump a number out of their tableaux, forming part of the matrix word $\pi$.
We avoid this by assuming in effect that our partitions have a previous row of
length $l=\delta_1$, which is empty in the skew shapes $\delta/\beta$ and
$\delta/\gamma$. 

\section{Internal Column Insertion, Hopscotch, and Stable Tableaux}\label{GNU}

For standard tableaux (and hence for all tableaux via standard renumbering),
internal column insertion is essentially the same as internal row
insertion---one simply replaces `row' by `column' in the definitions to
obtain internal column insertion.
While internal column insertion gives nothing new of itself, when combined
with complementation, we do get something interesting.
We call this new operation hopscotch.
Hopscotch is defined on pairs $P$ and $Q$, where $Q$ extends $P$, and one of
$P$, $Q$ is a tableau, while the other is a stable tableau, which is defined
in Section~\ref{sec:hopscotch} below.

\subsection{Local rules for column insertion}

We give the following
local rules formulation of internal column insertion for standard tableaux,
which is the matrix transpose of that for internal 
row insertion. 
Given partitions $\alpha\subdot\beta$ and $\alpha\subdot\gamma$, 
define a partition $\delta$ by
\begin{enumerate}
\item[(1)]
     If $\gamma\neq \beta$, then $\delta$ is their least upper bound,
     $\beta\vee\gamma$.
\item[(2)]
     If $\gamma=\beta$, and the box $\beta/\alpha$ is in the $j$th column,
     then  $\delta/\beta$ is a box in the ($j{+}1$)st column.
\end{enumerate}
Set
${\IC}(\alpha;\beta,\gamma):=\delta $.

This gives an algorithm for standard tableaux which can be generalised
to arbitrary tableaux by standard renumbering.  An explicit rule
for operating on horizontal strips is given in Section 5.2.

Consider now the reverse of this procedure.
Given partitions $\beta\subdot\delta$ and
$\gamma\subdot\delta$, we define $\alpha$ by
\begin{enumerate}
\item[(1)]
     If $\gamma\neq \beta$, then $\alpha$ is their greatest lower bound,
     $\beta\wedge\gamma$.
\item[(2)]
     If $\gamma=\beta$, and the box  $\delta/\beta$ is in the
     $j$th column, then $\beta/\alpha$ is a box in the ($j{-}1$)st column.
\end{enumerate}
Set
${\IC}(\beta ,\gamma ;\delta ):=\alpha $.

The reverse procedure does not work when $\delta/\beta$ is a box in the first
column.
This forces us to generalise our notions of shape and tableaux.  The
basic idea is to allow the creation of new columns, labeled with
non-positive integers, to the left of existing columns when we need them.  

Henceforth we define a {\it shape} to be a
finite weakly decreasing sequence of integers (positive or negative), called
{\it parts}.
This set of shapes with a fixed length $m$ forms a poset under componentwise
comparison, which extends Young's lattice of partitions.  Using $\subset$
for this partial order, we define skew shapes $\beta/\alpha$ for
$\alpha\subset\beta$ as before.
A horizontal strip is (as before), a skew shape $\beta/\alpha$ with 
at most one box in each column.  
We define a {\it tableau} of shape $\beta/\alpha$
with entries from $[k]$ to be a chain
$\alpha=\beta^0\subset\beta^1\subset\cdots\subset\beta^k=\beta$ where each
$\beta^i/\beta^{i{-}1}$ is a horizontal strip.  
We can convert a tableau $T$ (as defined here) into a tableau as
defined in~\cite{Sagan,Fu97,Stanley_ECII}, merely by adding a large enough
number to each part of 
the shapes that define $T$, in effect shifting $T$ horizontally.

\subsection{Internal column insertion}\label{sec:ici-def}
We extend this rule $\IC$ for standard tableaux to an insertion local rule
$\IC$ for tableaux so that the
resulting tableaux operation commutes with standard renumbering.
Let $\beta/\alpha$ and $\gamma/\alpha$ be horizontal strips, with
$\alpha$, $\beta$, and $\gamma$ shapes.
Consider applying the tableaux algorithm $\IC$ to the standard renumberings
of  the horizontal strips $\beta/\alpha$ and $\gamma/\alpha$.
This proceeds from left to right (from the last row to the first).
If there is any overlap between $\beta/\alpha$ and
$\gamma/\alpha$ in the $i$th row, then the entries in that row are shifted
to the right by the amount of that overlap, displacing entries in the
previous row, if necessary. 
For example, suppose $\alpha=(6,3,0)$, $\beta=(7,6,1)$, and
$\gamma=(8,4,1)$.
Then we have
$$
  \epsfbox{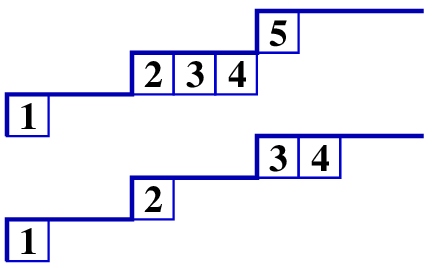}
   \qquad\raisebox{34pt}%
    {$\stackrel{\mbox{$\IC$}}%
      {\relbar\joinrel\relbar\joinrel\longrightarrow}$}
  \qquad
  \epsfbox{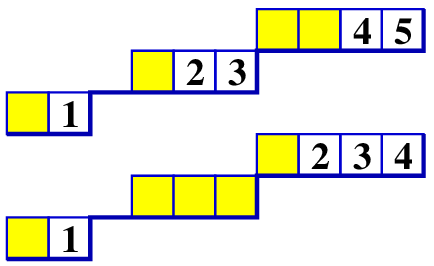} \ \ \raisebox{37pt}{.}
$$
and so $\delta=(10,6,2)$.

Formally, given shapes $\alpha$, $\beta$, and $\gamma$ with 
$\beta/\alpha$ and $\gamma/\alpha$ horizontal strips, we define the 
shape $\delta$ recursively as follows.
Suppose that the shapes $\alpha$, $\beta$, and $\gamma$ 
each have $m$ parts.
Set $d_m=0$, and for $i=m, m{-}1,\ldots,2$, set
\begin{equation}\label{IC}
  \begin{array}{rcl}
    \delta_i&:=&\min\left\{\begin{array}{l}
            \beta_i+\gamma_i+d_i-\alpha_i\\\alpha_{i{-}1}\end{array}\right.,\\
    d_{i{-}1}    &:=&\max\left\{\begin{array}{l}
            \beta_i+\gamma_i+d_i-\alpha_i-\alpha_{i{-}1}\\0\end{array}\right.,
  \rule{0pt}{22pt}\end{array}
\end{equation}
and set $\delta_1:=\beta_1+\gamma_1+d_1-\alpha_1$.
We define $\IC(\alpha;\beta,\gamma):=\delta$.
Note that $\IC$ does not change the number of rows of shapes.  
In the above example,  $d_3=d_2=0$ while $d_1=1$.
The numbers $d_i$'s keep track of boxes `bumped' from row $i{-}1$ to row $i$.

\begin{ex}\label{ex:ici}
We apply $\IC$ to the tableaux $P$ and $Q$ of
the upper left growth diagram in Figure~\ref{fig:four}
to obtain the following growth diagram.
\begin{equation}\label{eq:ici}
  \begin{array}{llll}
    1   & 2   & 31  & 33  \cr 
    21  & 31  & 42  & 53  \cr 
    211 & 311 & 421 & 531 \cr 
    221 & 321 & 431 & 541\cr 
  \end{array}
\end{equation}
\end{ex}

\begin{lem}
The rule~(\ref{IC}) defines a local rule $\IC$ which is reversible in that,
given $\beta,\gamma$, and $\delta$ with $\delta/\beta$ and $\delta/\gamma$
horizontal strips, there is a unique shape $\alpha$ with 
$\IC(\alpha;\beta,\gamma)=\delta$.  
\end{lem}

\noindent{\bf Proof. }
The only condition that needs to be checked is that $\delta$, $\beta$, and
$\gamma$ determine $\alpha$.
But this follows from the corresponding property of $\IR$ as we may
compute $\alpha$ 
from $\delta$, $\beta$, and $\gamma$ by the standard renumbering of
$\delta/\beta$ and $\delta/\gamma$, 
using the local rule ${\IC}$ for standard tableaux, which is equivalent to
$\IR$ via conjugation.
\QED

We define {\it internal column insertion} $\IC$ to be the tableaux operation 
determined by the local rule $\IC$.

\begin{thm}\label{thm:ici}
Internal column insertion is a bijection
$$
  \left\{
   \begin{minipage}[c]{1.7in}
    Tableaux $P$ and $Q$ which share an inner border.
   \end{minipage} 
  \right\}
  \ \stackrel{\IC}{\leftarrow\joinrel\relbar\joinrel\to}\ 
  \left\{
   \begin{minipage}[c]{1.7in}
    Tableaux $T$ and $U$ which share an outer border.
   \end{minipage} 
  \right\}
$$
with $\IC\circ\IC$ the identity 
such that if $\IC(P,Q)=(T,U)$, then  
$P$ is Knuth-equivalent to $T$ and $Q$ is Knuth-equivalent to  $U$.
Also, the outer border of $P$ equals the inner border of\/ $U$, and 
outer border of $Q$ equals the inner border of\/ $T$.

Furthermore, $\IC$-dual equivalence coincides with $\J$-dual
equivalence. 
\end{thm}

\noindent{\bf Proof. }
To show Knuth-equivalence,  let $\IC(P,Q)=(T,U)$.
Since $\IC$ commutes with standard renumbering, we may assume $Q$ is
standard. 
Then $T$ is obtained from $P$ by a series of internal column insertions.
Adapting the arguments of~\cite{SS90}, we see that internal column insertion
preserves Knuth equivalence.
Thus $P$ is Knuth-equivalent to $T$.
Since $\IC$ is symmetric in its two arguments, $Q$ is Knuth-equivalent
to $U$.  

Consider applying a sequence of $\IC$-moves to a tableau $P$.
Since $\IC$ commutes with standard renumbering, we may assume that each
$\IC$-move in that sequence is given by a tableau consisting of a single
box. 
The effect of these moves on the shape of $P$ will be the same as the effect
 on the shape of the standard renumbering $st(P)$ of $P$.
If we take the conjugate (matrix transpose) of tableaux in this
 sequence of $\IC$-moves applied to $st(P)$, we obtain a sequence of 
 $\IR$-moves applied to the conjugate $st(P)^t$ of $st(P)$.
Since $\J=\IR^C$, this observation and Theorem~\ref{R-de}
imply that $P$ and $Q$ are $\IC$-dual equivalent if and only if
 $st(P)^t$ and $st(Q)^t$ are $\J$-dual-equivalent.

However, a tableau $P$ is dual-equivalent to its standard renumbering, as
$\J$ commutes with standard renumbering.
Also, two dual-equivalent standard tableaux have dual equivalent conjugates,
as the conjugate of a jeu de taquin slide is another jeu de taquin slide,
by~(\ref{eq:js-def}).
This proves that $\IC $ and $\J $ have the same dual equivalence
classes.  
\QED

\subsection{Hopscotch and stable tableaux}\label{sec:hopscotch}

We would like to define the new tableaux operation $\HO$ of hopscotch to be
the tableaux operation complementary to internal column insertion $\IC$.
That is, if $P$ and $Q$ are tableaux with $Q$ extending $P$, then we would
set 
\begin{equation}\label{eq:ho-bad-def}
  \HO(P,Q)\ := (T,U^C)\qquad\mbox{where}\qquad
  \IC(P,Q^C)\ =\ (T, U)\,.
\end{equation}
(Recall that $U$ and $Q^C$ have the same content.)
Unfortunately, this does not work in general.
This is because $\IC$ can increase the number of columns in tableaux,
violating condition I of  Section~\ref{TaAl}.  
However, we will salvage something from this idea, by extending our
notion of ``tableau''.  

The essential problem comes from the fact that $\IC$ may increase
the number of columns of a tableau.
Thus, if we use a complementation parameter $l$ to form $Q^C$, 
then the tableau $U$ may extend beyond the $l$th column, and thus
forming $U^C$ will require a different complementation parameter, larger
than $l$. 
This cannot in general be remedied by increasing $l$ to some other
integer.
To solve this problem we build an asymmetry into $\HO$, requiring that
one of $P$ or $Q$ is a new object called a ``stable tableau'', which is
the complement of an ordinary tableau with respect to infinitely many
columns.  

We formalise this idea.
A {\it stable shape} $\alpha$ is a finite weakly decreasing sequence of
integers and the symbols $\infty$ and $\wbar{\infty}$,
where $\wbar{\infty}<n<\infty$ for any integer $n$.
If we regard a shape as having arbitrarily many trailing $\wbar{\infty}$s,
then stable shapes form a lattice under componentwise comparison, which
contains our earlier poset of shapes.
We form the stable skew shape $\beta/\alpha$ as before.
If the finite parts of $\alpha$ and $\beta$ occur in the same rows, then 
we can ignore initial $\infty$'s and trailing $\wbar{\infty}$'s and regard
$\beta/\alpha$ as an ordinary skew shape.
In this way, stable skew shapes and ordinary skew shapes
may extend each other.

A {\it stable horizontal strip} $\beta/\alpha$ is a pair of stable shapes
$\alpha\subseteq\beta$ where, if $\alpha_m$ is the first finite part of
$\alpha$ and $\beta_n$ the last finite part of $\beta$, then
$$
  \alpha_{m{-}1}=\infty=\beta_m>\alpha_m\geq\beta_{m{+}1}\geq\cdots\geq
  \alpha_{n{-}1}\geq\beta_n>\alpha_n=\wbar{\infty}=\beta_{n{+}1}\,.
$$
If we consider $\beta/\alpha$ as a collection of boxes in the plane, then
$\beta/\alpha$ has at most one box in each column, and it has 2 half-infinite
rows of boxes extending from either end.
Equivalently, if we let $\alpha^S:=(\infty,\alpha)$, then $\beta/\alpha$ is
a stable horizontal strip if and only if $\beta\subseteq\alpha^S$, and, upon
removing the $\infty$'s and $\wbar{\infty}$'s from $\beta$ and $\alpha^S$,
$\alpha^S/\beta$ is an ordinary horizontal strip with boxes in the empty
columns of $\beta/\alpha$.

\begin{figure}[htb]
$$
  \epsfxsize=164pt \epsfbox{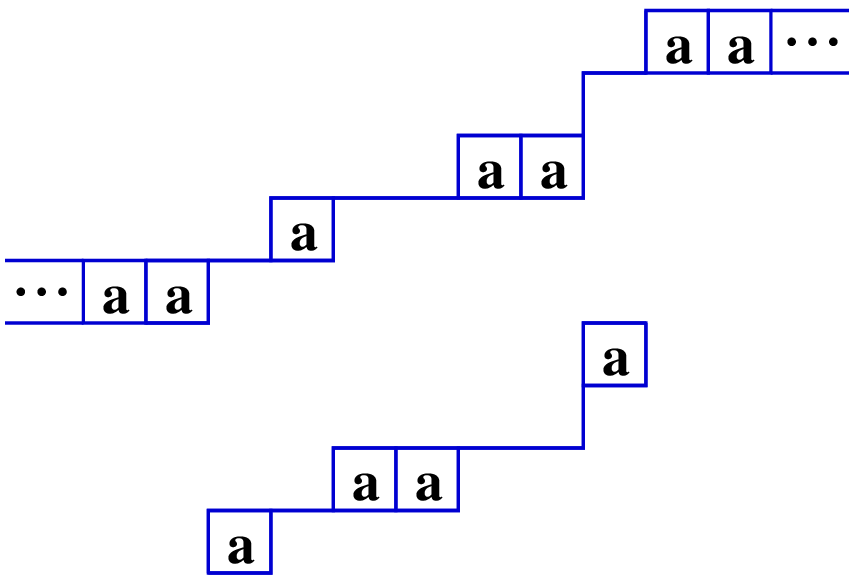}
$$
\caption{Complementary horizontal strips.\label{fig:shs}}
\end{figure}

Let $\alpha=(\infty,4,3,1,\bar{2},\wbar{\infty},\wbar{\infty})$ and 
$\beta=(\infty,\infty,3,3,\bar{1},\bar{3},\wbar{\infty})$, where we write
$\bar{a}$ for a negative integer $-a$.
Figure~\ref{fig:shs} shows  the stable horizontal strip $\beta/\alpha$ and
its complementary horizontal strip $\alpha^S/\beta$.

\begin{figure}[htb]
$$
  \epsfxsize=152pt \epsfbox{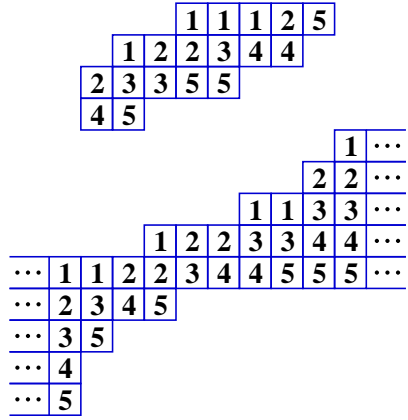}
$$
\caption{Stable complements.\label{fig:stable}}
\end{figure}

A {\it stable tableau} $T$ with stable shape $\beta/\alpha$ and entries from
$[k]$ is a  chain
$\alpha =\beta^{0}\subseteq \beta ^{1}\subseteq \cdots\subseteq
\beta^{k}=\beta$ of stable shapes,  where the
successive stable skew shapes $\beta^{i}/\beta^{i{-}1}$ are stable horizontal
strips.
The appropriate notion of complementation for stable tableaux involves
complementing all columns.
Thus if $T:\beta^{0}\subseteq \beta ^{1}\subseteq \cdots\subseteq
\beta^{k}$ is a stable tableaux, its {\it stable complement} $T^S$ is the chain
$\beta^k\subset(\infty,\beta^{k{-}1})\subset\cdots\subset(\infty^k,\beta^0)$.
Since $T$ is a stable tableau, each shape $\beta^i$ has one more part of
$\infty$ than its predecessor $\beta^{i{-}1}$ and the same number of finite
parts, and so we may remove the same number of parts of
$\infty$ from each shape $(\infty^{k{-}i},\beta^i)$ of $T^S$ and obtain an
ordinary tableau.
Similarly, if we complement an ordinary tableau $T$ using the
complementation parameter $\infty$, then we obtain its stable complement
$T^S$, a stable tableau.
Figure~\ref{fig:stable} shows a tableau $T$ and its stable
complement $T^S$ when $k=5$.

We summarise the properties of stable complementation.

\begin{thm}\label{thm:sc}
Let $T$ be a tableau of shape $\beta/\alpha$ and entries from $[k]$.
Then $T^S$ is a stable tableau with stable shape $(\infty^k,\alpha)/\beta$
and entries from $[k]$.
Similarly, if $T$ is a stable tableau with entries from $[k]$, then $T^S$ is
an ordinary tableau with entries from $[k]$.
In either case, we have $T^{SS}=T$.
\end{thm}

(For the last assertion, we identify tableaux that differ only by
a vertical shift.)

\begin{defn}\label{def:hop}
We define {\it hopscotch}, $\HO$, to be the stable complement of internal
column insertion.
That is, given a tableau $P$ and a stable tableau $Q$ with
either $P$ extending $Q$ or $Q$ extending $P$, then we set 
$$
  \HO(P;Q)\ := (T;U^S)\qquad\mbox{where}\qquad
  \IC(P,Q^S)\ =\ (T, U)\,.
$$

Suppose  we have
$\HO(P;Q)=(T;U)$.
We call the passage from the tableau $P$ to the tableau $T$ an 
{\it $\HO$-move
applied to $P$} (which depends upon the stable tableau $Q$).
Likewise, the passage from $Q$ to the stable tableau $U$ is also called 
a $\HO$-move.
As in Section~\ref{sec:dual}, this gives rise to the notion of $\HO$-dual
equivalence for tableaux and also for stable tableaux.
\end{defn}

\begin{thm}\label{thm:HO}
The tableaux operation hopscotch gives a bijection
$$
  \left\{
   \begin{minipage}[c]{1.9in}
    Pairs $(P;Q)$ with $P$ a tableau, $Q$ a stable tableau, and $Q$
   extending $P$.
   \end{minipage} 
  \right\}
  \ \stackrel{\HO}{\leftarrow\joinrel\relbar\joinrel\to}\ 
  \left\{
   \begin{minipage}[c]{1.9in}
     Pairs $(T;U)$ with $T$ a tableau, $U$ a stable tableau, and $T$
   extending $U$.
    \end{minipage} 
  \right\}
$$
with $\HO\circ\HO$ the identity.
If $\HO(P;Q)=(T;U)$, then $P$ is Knuth-equivalent to $T$.

Furthermore, two tableaux $P$ and $T$ are $\HO$-dual equivalent if and only
if $P$ and $T$ are $\J$-dual equivalent, and 
two stable tableaux $Q$ and $U$ are $\HO$-dual equivalent if and only
if the tableaux $Q^S$ and $U^S$ are $\J$-dual equivalent.
\end{thm}

\noindent{\bf Proof. }
These properties of hopscotch all follow from the corresponding properties
of internal column insertion (Theorem~\ref{thm:ici}), the definition of
hopscotch, and properties of stable complementation (Theorem~\ref{thm:sc}).
\QED

\section{Hopscotch and Tesler's Shift Games}

We give local rules
for hopscotch and relate it to Tesler's~\cite{Tesler}
rightward- and leftward-shift games.

\subsection{Local Rules for Hopscotch}
We give local rules for $\HO$ when the ordinary tableau is standard, that is
when each horizontal strip in $P$ is a single box.
As hopscotch commutes with standard renumbering (because internal column
insertion does), this enables the computation of $\HO(P;Q)$ when $P$ is an
arbitrary tableau.

\begin{defn}\label{def:loc-HO}
Let $\alpha\subset\beta\subdot\delta$  be shapes with 
$\beta/\alpha $ a stable horizontal strip.  Let $c$ denote the column of 
$\delta/\beta $, and $H$ denote the set of columns of 
$\beta/\alpha$.  
Define the stable shape $\gamma$ with 
$\alpha\subdot\gamma\subset\delta$ as follows.

\begin{enumerate}

\item[(1)]
     If $c\in H$, then define $\gamma $ so that $\alpha$ and $\gamma$ also
     differ in column $c$.  
     (So $\delta/\gamma$ will have the same set of columns $H$.)  

\item[(2)]
     If $c\not \in H$, then let $c'$ denote the smallest element of $H$
     greater than $c$, and define $\gamma$ so that $\gamma$ and $\alpha$
     differ in  column $c'$.  
     (So $\delta/\gamma $ will have the set of columns 
     $H':=H\cup \{c'\}\backslash \{c \}$.)
\end{enumerate}
\end{defn}

\begin{rem}\label{rem:HO-row}
In either case of Definition~\ref{def:loc-HO}, if $\beta$ and
$\delta$ differ in row $r$, then $\alpha$ and $\gamma$ will differ in row
$r'$, where $r'$ is the largest  row less than $r$ where $\alpha$ and
$\beta$ differ.
\end{rem}

\begin{figure}[htb]
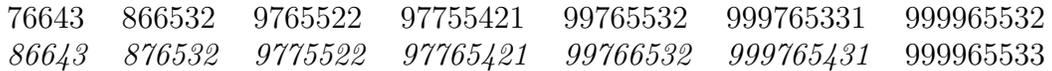

$$
  \begin{array}{lllllll}
 76643&\,866532&\,9765522&\,97755421&\,99765532&\,999765331&\,999965532 \cr 
\it 86643&\,\it 876532&\,\it 9775522&\,\it 97765421&
  \,\it 99766532&\,\it 999765431&\, 999965533
\end{array}
$$
\caption{Local rule of Definition~\ref{def:loc-HO}.\label{fig:2x6}}
\end{figure}

The local rule described in Definition~\ref{def:loc-HO} may sometimes be
computed when $\beta/\alpha$ is an ordinary horizontal strip---which may be
considered to be the truncation to a finite interval of columns of a stable
horizontal strip.
The $2\times 6$ array in Figure~\ref{fig:2x6} is a growth diagram of
ordinary shapes, with the second row filled in from right to left using the
local rule of Definition~\ref{def:loc-HO}.

\begin{thm}
Given stable shapes $\alpha\subset\beta\subdot\delta$ with $\beta/\alpha$ a
stable horizontal strip, the stable shape $\gamma$ defined in
Definition~\ref{def:loc-HO} satisfies
$\gamma=\HO(\alpha,\beta,\delta)$, i.e., Definition~\ref{def:loc-HO}
gives a local rules description of $\HO $.  
\end{thm}

\noindent{\bf Proof. }
The theorem follows from the description for $\IC$ when one strip
is standard.
Let $\alpha$, $\beta$, and $\gamma$ be shapes with $\alpha\subdot\gamma$,
$\alpha\subset\beta$ with $\beta/\alpha$ a horizontal strip.
Then the shape $\delta=\IC(\alpha;\beta,\gamma)$ may also be defined by
\begin{enumerate}
\item[(1)]
     If $\gamma\not\subset\beta$, then set $\delta:=\gamma\vee\beta$.
\item[(2)]
     Otherwise, $\gamma\subset\beta$.
     In this case, $\gamma/\alpha$ is a box in column $c$.
     Choose $\delta$ covering $\beta$ so that the box $\delta/\beta$ is in
     column $c'$ minimal subject to $c'>c$ with $\delta/\alpha$ is a
     horizontal strip.\QED
\end{enumerate}

The reverse local rule for $\HO$ given
$\alpha\subdot\beta\subset\delta$, when the inner horizontal strip consists
of a single box, is described by replacing Case 2 in
Definition~\ref{def:loc-HO} by 

\begin{enumerate}
\item[(2)$'$] \
     If $c\not \in H$, then let $c'$ denote the largest element of
     $H$ less than $c$, and define $\gamma$ so that $\gamma$ and $\alpha$
     differ in  column $c'$.  
     (So $\delta/\gamma $ will be the set of columns 
     $H':=H\cup \{c'\}\backslash \{c \}$.)
\end{enumerate}

\subsection{Hopscotch and Tesler's shift games}\label{sec:Tesler}
In his study of semi-primary lattices, Tesler~\cite{Tesler} defines certain
leftward- and rightward-shift games applied to standard skew tableaux which
model the effect of certain lattice-theoretic procedures.
These games give rise to an algorithm to construct a tableau $T$ of partition
shape from a standard skew tableau $P$.
Studying the corresponding objects in a semi-primary lattice, he then shows that
$T$ is the result of applying jeu de taquin slides to $P$, and thus $T$ is the
rectification of $P$.
We show how Tesler's algorithm may be regarded as a special case of
hopscotch and give a combinatorial (as opposed to geometric) proof that
his algorithm computes the rectification of $P$.
\smallskip

In these games of Tesler, a vertical
strip (full of $*$'s) is moved through 
a tableau $P$, and some entries of $P$ are removed (forming the $*$'s).
We describe the conjugate (replacing columns by rows)
of Tesler's rightward shift game
in terms of a local rule.
The leftward shift game is the reverse of this procedure.

\begin{defn}\label{def:rs}
An {\it almost standard tableau} $P$ is a tableau 
$\beta^0\subset\beta^1\subset\cdots\subset\beta^k$, where each 
$\subset$ in the chain is either a cover $\subdot$ or an equality.
Given an almost standard tableau
$P:\beta^0\subset\beta^1\subset\cdots\subset\beta^k$, Tesler's 
{\it rightward shift game} constructs another almost standard tableau
${\RS}(P):\alpha^0\subset\cdots\subset\alpha^k$ where each
$\beta^i/\alpha^i$ is a horizontal strip.
This begins with $\alpha^0:=\beta^0$.
If we have constructed $\alpha^0,\ldots,\alpha^{i{-}1}$, then we have the
(partial) growth diagram 
$$
  \begin{array}{cl}\alpha^{i{-}1}&\\\beta^{i{-}1}&\beta^i\end{array}\,.
$$

We construct $\alpha^i$ as follows.
\begin{enumerate}
\item[(1)]
     If $\beta^i=\beta^{i{-}1}$, then we set $\alpha^i:=\alpha^{i{-}1}$.
\item[(2)]
     If $\beta^i/\beta^{i{-}1}$ is a single box in the $r$th row and if the
     horizontal strip $\beta^{i{-}1}/\alpha^{i{-}1}$ has no boxes in rows less than
     $r$, then we set $\alpha^i:=\alpha^{i{-}1}$.
\item[(3)]
     If $\beta^i/\beta^{i{-}1}$ is a single box in the $r$th row and if the
     horizontal strip $\beta^{i{-}1}/\alpha^{i{-}1}$ has boxes in rows less
     than $r$, then we choose the largest row $r'$ of
     $\beta^{i{-}1}/\alpha^{i{-}1}$ less than $r$ and let
     $\alpha^i/\alpha^{i{-}1}$ be a box in that row. 
\end{enumerate}
\end{defn}

By Remark~\ref{rem:HO-row}, the similarity between these rules, particularly
(3), and those for hopscotch is evident. 

\begin{ex}\label{ex:RS}
We give a completed example of ${\RS}$ applied to a standard tableau
(the second row below):
$$
  \begin{array}{llllllllll}
   421 &\it  421 &\it  431 &\it  431 &\it  432 &\it  432 &\it  442 &\it  542 &\it  552 &\it  552 \\
   421 & 431 & 432 & 442 & 4421& 5421& 5422& 5522& 5532& 5542
  \end{array}
$$
We display this in terms of tableaux.
$$
   \epsfbox{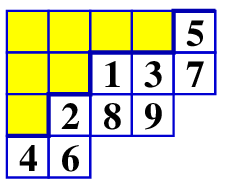}
   \quad\raisebox{24pt}%
    {$\stackrel{\mbox{${\RS}$}}%
      {\relbar\joinrel\relbar\joinrel\relbar\joinrel\longrightarrow}$}
  \quad
  \epsfbox{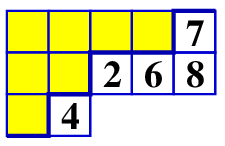}
$$
Informally, we move the entries in the tableaux in order northeast to
the nearest available space.  The 1 vacates the tableau, the 2 moves
into the empty space where the 1 was, the 3 vacates, the 4 moves where
the 2 was, the 5 vacates, the 6 moves where the 3 was, the 7 moves
where the 5 was, the 8 moves where the 7 was, and the 9 vacates.
\end{ex} 

Let $r(P)$ be the integers in $P$ but not in ${\RS}(P)$, that is,
the set of indices $i$ where $\beta^{i{-}1}\subdot\beta^i$ but 
$\alpha^{i{-}1}=\alpha^i$.
In the example above, $r(P)=\{1,3,5,9\}$.
Then Tesler's rectification algorithm runs as follows:
Given a standard skew tableau $P$, form tableaux
$P_0,\ldots,P_k$ and sets $r_1,\ldots,r_k$ recursively, 
initially setting $P_0:=P$ and then 
$$
  P_{i}\ :=\ {\RS}(P_{i{-}1})\qquad\mbox{and}\qquad
  r_i\ :=\ r(P_{i{-}1})\,.
$$

\begin{prop}[\cite{Tesler}, Theorem 8.13]\label{prop:tesler}
The tableau $P_{k}$ is  empty and $r_1,\ldots,r_k$ form the
rows of the rectification $S$ of $P$.
\end{prop}

We relate this to hopscotch, first showing that ${\RS}$ is
equivalent to a particular $\HO$-move applied to 
$P$, and then we show how to compute the tableau $S$ of
Proposition~\ref{prop:tesler} using hopscotch.
This will give a new, combinatorial proof of this result of Tesler.

Let $P:\beta^0\subdot\beta^1\subdot\cdots\subdot\beta^k$ be a standard
tableau where each shape $\beta^i$ has $n$ nonnegative parts (appending
0's if needed).  
Let $b$ be the first (largest) part of $\beta^k$.
Prepend $\infty$ and append 0 to each shape $\beta^i$ and consider it to
be a stable shape.
Set $\alpha^0=\beta^0$ and define
$$
  \gamma^0\ :=\ b\geq b_1\geq b_2\geq \cdots\geq b_n>\wbar{\infty}\,,
$$
where $\beta^0=\infty>b_1\geq\cdots\geq b_n\geq 0$.
Then the stable horizontal strip $\beta^0/\gamma^0$ consists of two
semi-infinite rows of boxes, one beginning in column $b{+}1$ in the first row,
and one ending in column 0 in row $n{+}1$.
Set $Q:=\beta^0/\gamma^0$ and $(T;U):=\HO(P;Q)$. 
Define shapes $\gamma^1,\ldots,\gamma^k$ and $\alpha^0,\ldots,\alpha^k$ so
that the standard tableau $T$ above is
$\gamma^0\subdot\gamma^1\subdot\cdots\subdot\gamma^k$ and 
${\RS}(P)=\alpha^0\subset\cdots\subset\alpha^k$.
For each $i=0,1,\ldots,k$, set 
$a_i:=\#\{j<i\mid \alpha^j=\alpha^{j{+}1}\}$.

\begin{lem}\label{lem:tesler}
For each $i=0,\ldots,k$ the stable shape $\gamma^i$ is 
\begin{equation}\label{eq:tes}
   b+a_i\geq\alpha^i_1\geq\cdots\geq\alpha^i_n>\wbar{\infty}\,.
\end{equation}
In particular, $r(P)$ is the first row of\/ $T$ and ${\RS}(P)$
consists of the remaining rows of\/ $T$.
\end{lem}

\noindent{\bf Proof. }
We prove this by induction on $i$, the case $i=0$ being the definition of
$\gamma^0$.
Suppose that~(\ref{eq:tes}) holds for $\gamma^0,\ldots,\gamma^{i{-}1}$.
We compare the $i$th steps of $\HO$ and ${\RS}$.
Since $P$ is standard, case (1) of Definition~\ref{def:rs} for
${\RS}$ does not occur. 
If we are in case (2) of Definition~\ref{def:rs}, then
$\alpha^i=\alpha^{i{-}1}$ and $a_i=1+a_{i{-}1}$.
Since the stable horizontal strip $\beta^{i{-}1}/\gamma^{i{-}1}$ has no
boxes in rows between its first row and the row $r$ of the single box 
$\beta^i/\beta^{i{-}1}$,  the stable shapes $\gamma^i$ and $\gamma^{i{-}1}$
differ only in the first row, by Remark~\ref{rem:HO-row}.
Finally, case (3) of Definition~\ref{def:rs} is equivalent to hopscotch.
\QED

\begin{figure}[htb]
$$
  \epsfbox{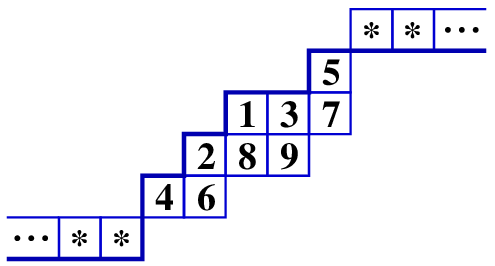}
   \quad\raisebox{34pt}%
    {$\stackrel{\mbox{$\HO$}}%
      {\relbar\joinrel\relbar\joinrel\longrightarrow}$}
  \quad
  \epsfbox{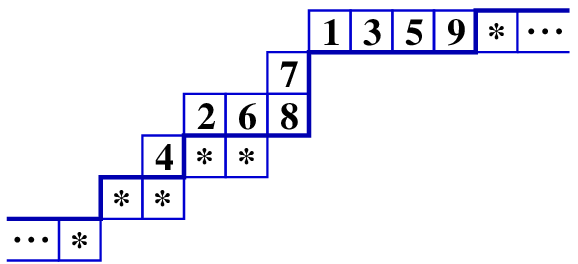}
$$
\caption{${\RS}$ as a hopscotch move.\label{fig:tes-row}}
\end{figure}

Figure~\ref{fig:tes-row} shows an example of this hopscotch move where
we fill the stable horizontal strip $Q$ with  $*$'s. 
Note the similarity with Example~\ref{ex:RS}.

\begin{figure}[htb]

$$
  \epsfbox{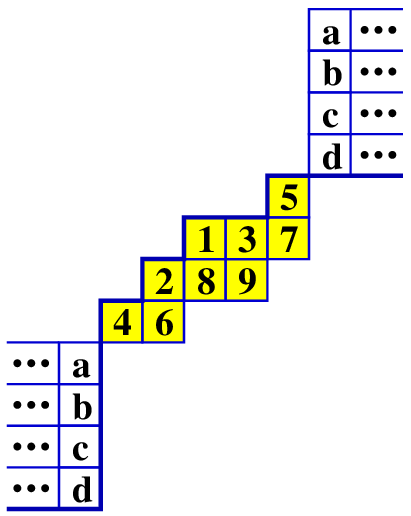}
   \qquad\raisebox{55pt}%
    {$\stackrel{\mbox{$\HO$}}%
      {\relbar\joinrel\relbar\joinrel\longrightarrow}$}
  \qquad
  \epsfbox{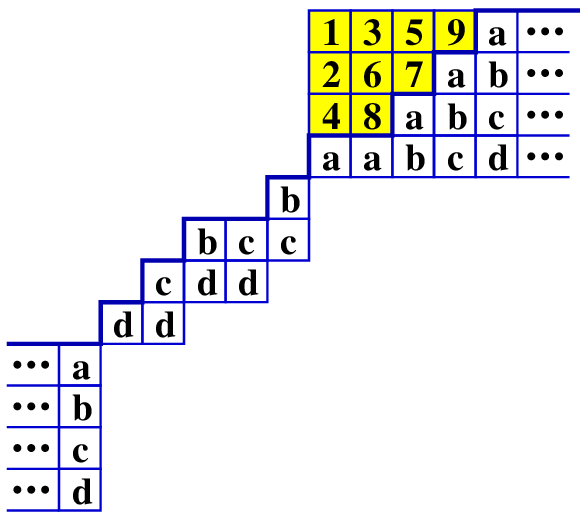}
$$
\caption{Tesler's rectification algorithm as hopscotch.\label{fig:HO-RS}}
\end{figure}

To compute the rectification tableau $S$ of Proposition~\ref{prop:tesler}
using hopscotch, we iterate the hopscotch move of Lemma~\ref{lem:tesler},
letting the initial tableau $P$ extend sufficiently many stable horizontal
strips of the form $Q$ above.
Given $P$ with $\beta^i$, $b$, and $b_i$ as in the paragraph before
Lemma~\ref{lem:tesler}, set $l:=\min\{k,n\}$.
Define the stable tableau
$Q:=\alpha^0\subset\alpha^1\subset\cdots\subset\alpha^l$ by
$$
  \alpha^i\ :=\ \infty^i>b^{l{-}i}\geq b_1\geq\cdots\geq b_n\geq 0^i
                > \wbar{\infty}^{l{-}i}\,.
$$
If we prepend $\infty^l$ and append $0^l$ to each partition $\beta^i$, then
the standard tableau $P$ extends the stable tableau $Q$.
Set $(T;U):=\HO(P;Q)$. 
Figure~\ref{fig:HO-RS} displays this action of
hopscotch on the tableau $P$ of Figure~\ref{fig:tes-row}, rectifying it
to $T$.  
The tableau $T$ has partition shape and is in the initial $l$ rows and
columns greater than $b_n$.
By Theorem~\ref{thm:HO}, $T$ is Knuth-equivalent to $P$ and hence is the
rectification of $P$.
Successively applying Lemma~\ref{lem:tesler} and analysing the effect of
hopscotch on the columns greater than $b_n$, we obtain a combinatorial
proof that Tesler's algorithm rectifies tableaux. 

\section{Jeux de Tableaux}
In Section~\ref{Tableaux} we asserted that the new tableaux operations of
$\IR$, $\IC$, and $\HO$ have ramifications for ordinary tableaux, even
though we generalised ordinary tableaux in order to define these
operations. 
We describe four different algorithms to rectify a column-strict tableau,
supplementing the classical jeu de taquin of Sch\"utzenberger~\cite{Schu77}.
The first is constructed from the internal (row) insertion of Sagan and
Stanley~\cite{SS90}, and the second is similarly derived from column
insertion~\cite{Fu97,Sagan}. 
The last two are related to hopscotch.
One, which we call row extraction, is the column-strict tableaux version of
Tesler's rectification algorithm of Section 6.2, and the other, which we
call column extraction, is adapted from Stroomer's column sliding
algorithm~\cite{Stroomer}. 

In this section, tableaux all have (skew) partition shape, in the
traditional sense: non-negative integer parts, and the initial part of a
partition $\lambda$ is $\lambda_1$.

\subsection{Internal Insertion Games}
Let $P$ be a column-strict tableau of shape $\lambda/\mu$ with
$\mu,\lambda$ partitions.
We assume that $P$ has entries in both the first row and first
column.
A(n)  {\it (inside) cocorner} is an entry of $P$ with no neighbours in $P$
to the left or above.
An {\it internal row insertion} on $P$ as introduced in~\cite{SS90}
begins with a cocorner of $P$.
That entry is removed from $P$ and inserted into the subsequent rows of $P$
using Schensted (row) insertion.
The resulting tableau is Knuth-equivalent to $P$.

\begin{figure}[htb]
$$
  \epsfbox{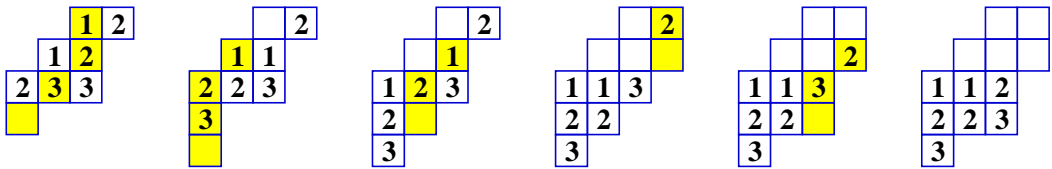}
$$
\caption{The row insertion game.\label{fig:rectIR}}
\end{figure}

The {\it row insertion game} begins with a column-strict tableau $P$ of
shape $\lambda/\mu$ and a row $j$ with $\mu_j=0$.
The game proceeds by successive internal row insertions beginning with
cocorners in rows $<j$, and it ends when there are no such cocorners.
The result is a tableau $T$ of partition shape (occupying rows $\geq j$)
Knuth-equivalent to $P$, that is, the rectification of $P$.
Thus the result of the row insertion game is independent of
the particular sequence of cocorners chosen.
We illustrate this process in Figure~\ref{fig:rectIR}.
In each tableau in that sequence, we shade the cocorner and insertion path
that creates the subsequent tableau.\smallskip

We do the same with column insertion.
An {\it internal column insertion} on a tableau $P$ begins with a cocorner
of $P$.
That entry is removed from $P$ and inserted into subsequent columns of $P$
using column insertion~\cite{Fu97,Sagan}, and 
the resulting tableau is Knuth-equivalent to $P$.

\begin{figure}[htb]
$$
  \epsfxsize=302.67pt\epsfbox{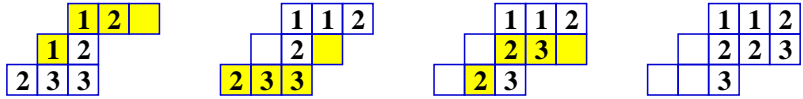}
$$
\caption{The column insertion game.\label{fig:rectIC}}
\end{figure}

The {\it column insertion game} begins with a column-strict tableau $P$ of
shape $\lambda/\mu$
and proceeds by successive internal column insertions beginning with
cocorners in columns $1,\ldots,\mu_1$, and it ends when there are no such
cocorners. 
The result is a tableau $T$ of partition shape (occupying columns
$>\mu_1$) 
Knuth-equivalent to $P$, that is, the rectification of $P$.
Thus the result of the column insertion game is independent of
the particular sequence of cocorners chosen.
We illustrate this process in Figure~\ref{fig:rectIC}.
In each tableau in that sequence, we shade the cocorner and insertion path
that creates the subsequent tableau.

\subsection{Row Extraction}
Row extraction is the extension of the rectification algorithm of
Tesler, ${\RS}$, 
described in Section~\ref{sec:Tesler}, to column-strict tableaux.
After describing this extension, we show how it may be
used to compute the rectification of a column-strict skew tableau $P$.
Define the {\it cross order} on the cells of a shape or tableau by
$c\prec d$ if the cell $c$ lies in the same column as, or a subsequent
column to, $d$  and is in a previous row.

Let $P: \mu=\lambda^0\subset\lambda^1\subset\cdots\subset\lambda^k=\lambda$ be a 
tableau where $\lambda^i/\lambda^{i{-}1}$ is a horizontal strip of $i$'s.
{\it Row extraction} creates an initially vacant horizontal strip of $*$'s and
moves it successively past each horizontal strip $\lambda^i/\lambda^{i{-}1}$ of
$P$, possibly getting larger as it goes by the creation of new $*$'s, which  
replace entries of $P$.
The result is a tableau ${\RE}(P)$ with inner
border $\mu$ extended by a horizontal strip of $*$'s whose outer border
is $\lambda$.
The collection of all replaced entries of $P$ is a multiset $r(P)$.

\begin{figure}[htb]
$$
  \epsfxsize=110.33pt\epsfbox{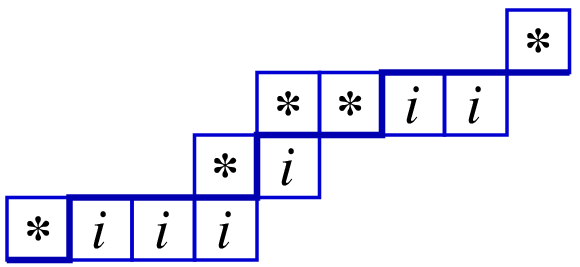}
   \qquad\raisebox{12pt}%
    {$\relbar\joinrel\relbar\joinrel\longrightarrow$}
  \qquad
  \epsfxsize=110.33pt\epsfbox{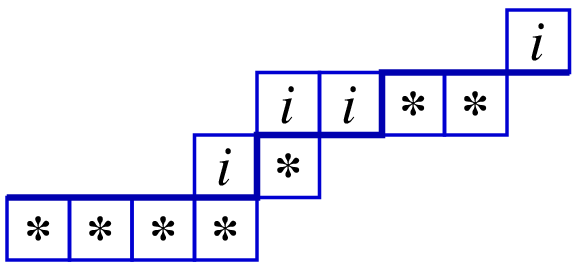}
$$
\caption{Row extraction on horizontal strips.\label{fig:RE}}
\end{figure}

We describe how to move a horizontal strip of $*$'s past a horizontal strip 
$\lambda^i/\lambda^{i{-}1}$ extending it.
This proceeds left-to-right through $\lambda^i/\lambda^{i{-}1}$, moving each entry 
$i$ as follows. 
The current entry $i$ in $\lambda^i/\lambda^{i{-}1}$ is interchanged with the
maximal $*$ smaller than it in the cross order, and if there is no such $*$,
then that $i$ is removed and replaced by a new $*$.
The horizontal strip of $*$'s is initially empty, and remains so until 
encountering the first nonempty horizontal strip $\lambda^i/\lambda^{i{-}1}$, 
at which point all $i$'s are logically removed and replaced by $*$'s.
We give an example in Figure~\ref{fig:RE} of an intermediate stage of the algorithm, where two $i$'s are removed.

\begin{figure}[htb]
$$
  \begin{array}{l}
     \epsfbox{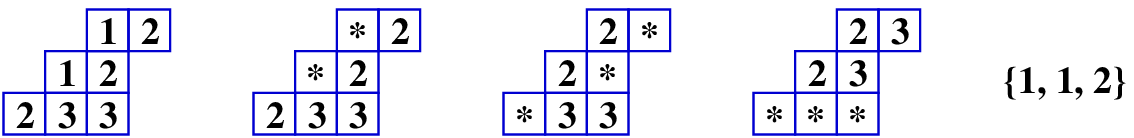}\\
     \epsfbox{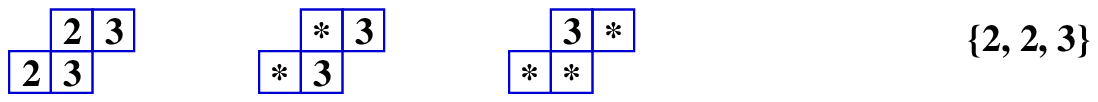}\\
     \epsfbox{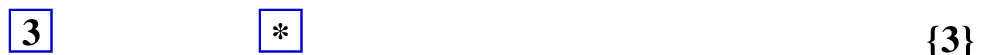}
   \end{array}
$$
\caption{Rectifying using row extraction.\label{fig:rectRE}}
\end{figure}

This process is essentially Tesler's algorithm ${\RS}$ as described in
Section~\ref{sec:Tesler}, up to standard renumbering.
The column-strict extension of Tesler's rectification algorithm proceeds as follows:
Given $P$, form tableaux $P_0$, $P_1$, $\ldots$, $P_k$ and
multisets $r_1,\ldots,r_k$ by $P_0:=P$, and for each $i=1,\ldots,k$, 
$$
   P_i\ :=\ {\RE}(P_{i{-}1})\qquad\mbox{and}\qquad 
   r_i\ :=\ r(P_{i{-}1})\,.
$$
Then $P_k=\emptyset$, and an analysis as in Section~\ref{sec:Tesler}
shows that the  
multisets $r_1,\ldots,r_k$ form the rows of the rectification of $P$.
Figure~\ref{fig:rectRE} shows this procedure applied to the tableau $P$ of
Figures~\ref{fig:rectIR} and~\ref{fig:rectIC}.
Each row in the figure is one application of ${\RE}$, and the multiset of entries
removed  in each step is displayed to the right.

\subsection{Column Extraction}
The last algorithm of column extraction is strangest of all.
Its discovery in February 1992 and a  
desire to find a common framework with jeu de taquin was the
genesis of this paper.
Column extraction is constructed from the column sliding algorithm of
Stroomer~\cite{Stroomer}, which we first describe.

Let $P$ be a column-strict tableau with entries in $[k]$, and let $b$ be an
inner corner of $P$.
A {\it column slide} beginning at $b$ is given by moving an empty box $b$
through $P$ as follows.
Initially $b$ switches with the first 1 greater in the cross order in
$P$. 
Then the empty box switches with the first 2 in the cross order, and so
on, concluding when $b$ switches with a $k$. 
It is only permissible to begin a column slide at $b$ when there exists
a chain of entries in the cross order labelled $1,2,\ldots,k$.  
In the figure below, it is not permissible to begin a column slide at the
inner corner marked with an $X$, but it
is permissible to begin a column 
slide at the inner corner marked with an $*$.
We shade the paths of both the impermissible column slide beginning with
$X$ and the permissible column slide beginning with $*$, and then in the
right tableau show the result 
of the permissible column slide.
\begin{figure}[htb]
$$
  \epsfbox{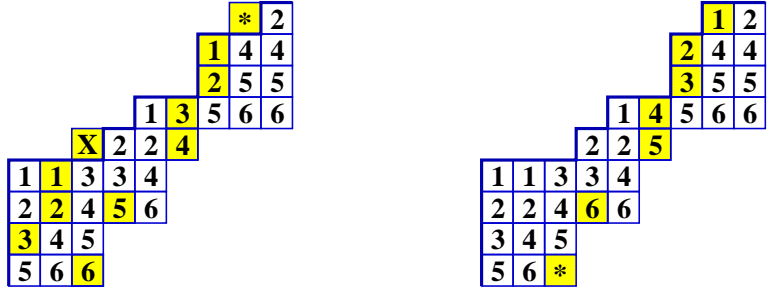} 
$$
\caption{Impermissible and  permissible column slides.\label{fig:CS}}
\end{figure}

\begin{lem}\label{lem:CS}
The result $T$ of a permissible column slide on a tableau $P$ is 
Knuth-equivalent to $P$.
\end{lem}

\noindent{\bf Proof. }
Consider the complement $P^C$  writing $P^C$ above the
tableau $P$.
Then an inner corner $b$ of $P$ is an outer cocorner of $P^C$ from which we
could begin a reverse column insertion~\cite{Fu97,Sagan}.
Either an entry of $P^C$ will be bumped from $P^C$ or this reverse insertion
will result in a new box on the inner border of $P^C$.
Stroomer~\cite{Stroomer} showed that in the first case, it is not
permissible to begin a column slide at $b$ on $P$, and in the second case,
the column slide is permissible, and $T^C$ is the tableau obtained from the
(internal) column insertion.
Thus $P^C$ and $T^C$ are Knuth-equivalent, and so by
Theorem~\ref{comp-equiv}, $P$ and $T$ are Knuth-equivalent.
\QED

This proof shows that a permissible column slide is complementary to an
internal column insertion, and thus may be computed using hopscotch moves.
Indeed, the left and right tableaux in Figure~\ref{fig:CS} are respectively
the second and first rows in Figure~\ref{fig:2x6}.
We illustrate the complementary nature of a permissible column slide and an
internal insertion. 
The (reverse) column insertion in the tableau below on the left
beginning with the shaded box results in the tableau on the right below.
These tableaux are the complements of the tableaux in Figure~\ref{fig:CS}.
$$
  \epsfbox{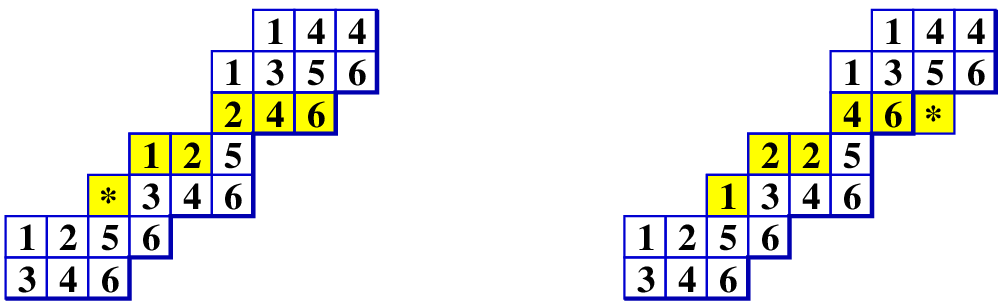}
$$

Observe that if the first column of the tableau $P$ is full, that is,
consists of all the entries in $[k]$, then any inner corner of $P$ can
initiate a permissible column slide.
Also note that if $T$ is the result of a permissible column slide beginning
at an inner corner $b$, then we may initiate a permissible column slide 
at any inner corner of $T$ smaller than $b$ in the
cross order.

\begin{defn}\label{def:CE}
Let $P$ be a tableau with entries in $[k]$.
Form the tableau $P'$ by placing a full column of the entries $1,\ldots,k$
to the left of, and beginning in the same row as, the first column of $P$.
The cell $b$  just above the rightmost column of $P'$ is an inner corner of
$P'$, and so it is permissible to initiate a column slide in that cell.
The cell just above $b$ is an inner corner of the resulting tableau, and we may begin another column slide with this cell.
Repeating this procedure at most $k$ times, we obtain a tableau $T'$ whose
rightmost column is full, and we delete this column to obtain the tableau
${\CE}(P)$. 
\end{defn}

By construction, $P$ and ${\CE}(P)$ have the same content.
More is true.

\begin{lem}\label{lem:CE-KE}
The tableaux $P$ and  ${\CE}(P)$ are Knuth-equivalent.
\end{lem}

We will prove this lemma at the end of this section.
We illustrate this operation on the tableau of our running example.
\begin{figure}[htb]
$$
  \epsfbox{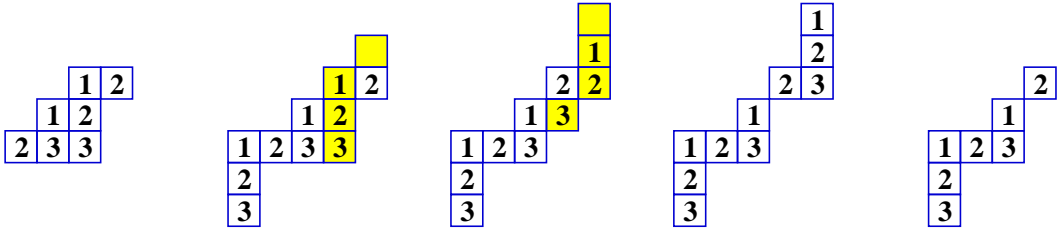}
$$\caption{Column extraction.\label{fig:CE}}
\end{figure}

Iterating ${\CE}$ sufficiently many times rectifies a tableau.

\begin{thm}\label{thm:rectCE}
Let $P$ be a column-strict tableau.
Let $j$ be the number of columns of $P$ that begin in rows higher than the
first column of $P$.
Then the $j$th iterate ${\CE}^j(P)$ of column extraction applied to
$P$ is the rectification of $P$.
\end{thm}

\noindent{\bf Proof. }
The number of columns of ${\CE}(P)$ is at most the number
of columns of $P$.
Also, if the first $i$ columns of $P$ begin in the same row, then the first
$i{+}1$ columns of ${\CE}(P)$ will begin in that same row, unless 
${\CE}(P)$ has fewer than $i{+}1$ columns.
Thus ${\CE}^j(P)$ has all of its columns beginning in the same row,
and so it has partition shape.
Since $P$ is Knuth-equivalent to ${\CE}^j(P)$ by
Lemma~\ref{lem:CE-KE}, we see that ${\CE}^j(P)$ is the
rectification of $P$.
\QED

\begin{figure}[htb]
$$
   \begin{array}{l}
    \epsfbox{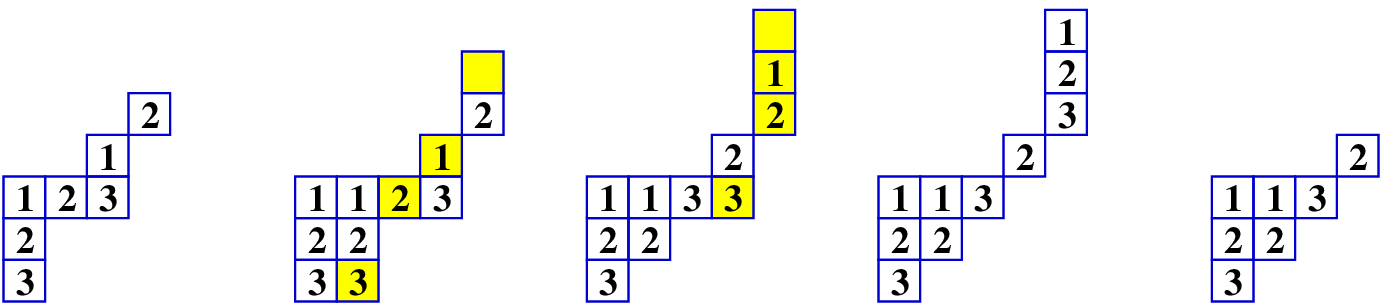}\\
    \epsfbox{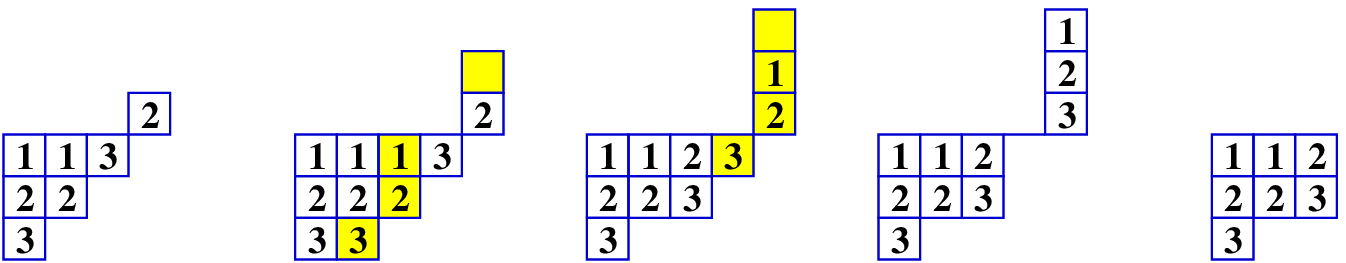}
   \end{array}
$$
\caption{Rectifying using column extraction.\label{fig:rectCE}}
\end{figure}

Figure~\ref{fig:rectCE} illustrates the algorithm of
Theorem~\ref{thm:rectCE}, applying ${\CE}$ twice more 
to the example of Figure~\ref{fig:CE}.

Our proof of Lemma~\ref{lem:CE-KE} uses Knuth equivalence $\simeq$ of words
in the alphabet $[k]$~\cite{Sagan,Fu97,Stanley_ECII}.
Given two words $w,v$, let $w.v$ be their concatenation.
Given a column-strict tableau $P$, let $\omega(P)$ be its 
{\it (column) word}, which is the entries of $P$ listed from the bottom to
top in each column, starting in the leftmost column and moving right.
We use the result of Sch\"utzenberger~\cite{Schu77} that tableaux $P$
and $Q$ are Knuth-equivalent if and only if $\omega(P)\simeq\omega(Q)$.
Set $C:=k\ldots 2.1$, the word of a full column.
We first prove a lemma concerning Knuth equivalence and $C$, or rather
commutation and cancellation in the plactic monoid of Lascoux and
Sch\"utzenberger~\cite{LS81}. 

\begin{lem}\label{lem:col-KE}
Let $w,v$ be any words in the alphabet $[k]$.
Then
\begin{enumerate}
\item[(\rm i)]
      $C.w\simeq w.C$.
\item[(\rm ii)]
      $C.w\simeq C.v$ if and only if $w\simeq v$.
\end{enumerate}
\end{lem}

\noindent{\bf Proof. }
For (i), it suffices to consider the case when $w$ is a single number $i$.
Then $i.C$ is the column word of a tableau of skew shape 
$2^k/1^{k{-}1}$ whose rectification has partition shape $(2,1^{k{-}1})$ and word 
$C.i$.

The reverse direction of (ii) is trivial.
Suppose $C.w\simeq C.v$ and let $P,Q$ be tableaux of partition shape with
$\omega(P)\simeq w$ and $\omega(Q)\simeq v$.
Then $C.\omega(P)\simeq C.\omega(Q)$.
But for any tableau $T$, $C.\omega(T)$ is the word of a tableau $T'$ of
partition shape obtained from $T$ by placing a full column to its left.
Since there is a unique tableau of partition shape in any Knuth equivalence
class, we must then have $C.\omega(P)= C.\omega(Q)$, thus $P'=Q'$,  and so
$P=Q$. 
But this implies $w\simeq v$.
\QED

\noindent{\bf Proof of Lemma~\ref{lem:CE-KE}. }
Let $P,P',T'$ be as in Definition~\ref{def:CE} and set  
$T={\CE}(P)$.
Then 
$$
   C.\omega(P)\ =\ \omega(P')\ \simeq\ \omega(T')\ 
   =\ \omega(T).C\ \simeq C.\omega(T)\,.
$$
The first Knuth-equivalence follows from Lemma~\ref{lem:CS}, as $T'$ is
obtained from $P'$ by column slides.
The second Knuth-equivalence is Statement (i) of Lemma~\ref{lem:col-KE}.
Since $C.\omega(P)\simeq C.\omega(T)$, we deduce that 
$\omega(P)\simeq \omega(T)$, by (ii) of Lemma~\ref{lem:col-KE}, and so $P$ is Knuth equivalent to 
${\CE}(P)$.
\QED

\end{document}